\documentclass[10 pt, journal, final, twocolumn]{IEEEtran}

\usepackage{enumerate}
\usepackage{amsmath, amsthm, amssymb}
\usepackage{float}
\usepackage{subcaption}
\usepackage{flexisym}
\usepackage{mathtools}
\usepackage{breqn}
\usepackage{array}
\usepackage{tikz}
\usepackage{cite}
\usetikzlibrary{shapes.geometric, arrows}
\usepackage[font = small]{caption}
\usepackage{paralist}
\usepackage{mathrsfs}
\usepackage{algorithm}
\usepackage{algpseudocode}
\usepackage[inline]{enumitem}
\allowdisplaybreaks[1]

\newtheorem{proposition}{Proposition}
\newtheorem{lemma}{Lemma}
\newtheorem{theorem}{Theorem}
\newtheorem{definition}{Definition}

\newtheorem{assumption}{Assumption}
\newtheorem{specification}{Specification}

\DeclareMathOperator*{\trace}{trace}
\DeclareMathOperator*{\diag}{diag}
\DeclareMathOperator*{\blkdiag}{blkdiag}

\renewcommand{\IEEEQED}{\IEEEQEDclosed}

\begin{document}
\title{Resilient Distributed Estimation Through Adversary Detection}

\author{Yuan Chen, Soummya Kar, and Jos\'{e} M. F. Moura
\thanks{Yuan Chen {\{(412)-268-7103\}}, Soummya Kar {\{(412)-268-8962\}}, and Jos\'{e} M.F. Moura {\{(412)-268-6341, fax: (412)-268-3890\}} are with the Department of Electrical and Computer Engineering, Carnegie Mellon University, Pittsburgh, PA 15217 {\tt\small \{yuanche1, soummyak, moura\}@andrew.cmu.edu}}
\thanks{This material is based upon work supported by the Department of Energy under Award Number DE-OE0000779 and by DARPA under agreement numbers DARPA FA8750-12-2-0291 and DARPA HR00111320007. The U.S. Government is authorized to reproduce and distribute reprints for Governmental purposes notwithstanding any copyright notation thereon. The views and conclusions contained herein are those of the authors and should not be interpreted as necessarily representing the official policies or endorsements, either expressed or implied, of DARPA or the U.S. Government.}}
\maketitle

\begin{abstract}
	This paper studies resilient multi-agent distributed estimation of an unknown vector parameter when a subset of the agents is adversarial. We present and analyze a Flag Raising Distributed Estimator ($\mathcal{FRDE}$) that allows the agents under attack to perform accurate parameter estimation and detect the adversarial agents. The $\mathcal{FRDE}$ algorithm is a consensus+innovations estimator in which agents combine estimates of neighboring agents (consensus) with local sensing information (innovations). We establish that, under $\mathcal{FRDE}$, either the uncompromised agents' estimates are almost surely consistent or the uncompromised agents detect compromised agents if and only if the network of uncompromised agents is connected and globally observable. Numerical examples illustrate the performance of $\mathcal{FRDE}$.
\end{abstract}


\section{Introduction}\label{sect: intro}
Distributed algorithms arise in numerous applications such as consensus algorithms~\cite{TBA, Boyd1}, inference and computation over wireless networks~\cite{Kar1, Kar2, Kar3, Sundaram2, DistributedWSN, DiffusionLMS}, state estimation in the electric power grid~\cite{Xie2}, and control of multi-agent systems~\cite{CooperationSurvey, InterconnectedSystems}. In these applications, individual agents exchange information with their neighbors to compute the average of a set of initial agent values~\cite{TBA, Boyd1} or to rendezvous at a common location~\cite{InterconnectedSystems}. 

Distributed algorithms are vulnerable to attack from adversarial agents. The presence of compromised agents can prevent a distributed algorithm from achieving its desired goal. For example, one popular algorithm for distributed parameter estimation is the Diffusion Least Mean Squares (DLMS) algorithm~\cite{DiffusionLMS}. In the DLMS algorithm, like in other general distributed estimation methods, the presence of adversarial agents can prevent the other agents from correctly estimating the parameter of interest~\cite{LMSCluster1}. References~\cite{LMSCluster1} and~\cite{LMSCluster2} propose a distributed technique, coupled with diffusion adaptation, to detect and counteract adversarial agents who behave as though they observe a different parameter than the true parameter of interest. The modifications presented in~\cite{LMSCluster1} and~\cite{LMSCluster2} address a specific type of adversarial behavior, but, in general, there are many other types of attacks. 

This paper addresses the following question: how can normally behaving agents estimate an unknown parameter in the presence of adversarial agents? Consider a network of agents, each making local sensor measurements of an unknown parameter $\theta^*$. While the agents have enough information, collectively, to determine $\theta^*$ (i.e., agent models are globally observable), no individual agent has enough information by itself to determine $\theta^*$ (i.e., agent models are locally unobservable).\footnote{In contrast, DLMS (\cite{DiffusionLMS, LMSCluster1, LMSCluster2}) require agent models to be locally observable.} The agents exchange information with their direct neighbors to perform estimation. In this paper, we propose a method for agents to either detect adversaries or achieve resilient distributed estimation. In other words, we establish a necessary and sufficient condition such that, either, the attack is weak to escape detection (with false alarm probability below a desired level), in which case our algorithm is resilient and the estimates converge almost surely to $\theta^*$, or, the attack is strong, in which it is detected.

\subsection{Literature Review}
An important example that demonstrates the effect of misbehaving agents in distributed computation is the Byzantine Generals Problem~\cite{ByzantineGenerals} where a group of agents must decide, by passing messages between one another (in an all-to-all manner), whether or not to attack an enemy city. Adversarial agents attempt to disrupt the message passing process and cause the group to reach the incorrect decision. The authors of~\cite{ByzantineGenerals} show that, if the number of adversarial agents is at least one third of the number of all agents, then it is impossible to design any distributed algorithm (under the specified admissible class of message passing protocols) to compute the correct decision. Conversely, if the number of adversarial agents is less than one third of the number of all agents, they provide an algorithm for the non-adversarial agents to reach the correct decision even in the presence of adversaries. Reference~\cite{ByzantineGenerals} addresses consensus (i.e., agreeing on a decision whether or not to attack) in an all-to-all communication setting. In contrast, this paper studies distributed estimation in a sparse communication setting (i.e., agents may only communicate with their neighbors).

Existing work has addressed Byzantine attacks in the context of decentralized inference~\cite{Varshney1, Varshney2, FusionCenterEstimation1, FusionCenterDetection2}, where individual agents make measurements of an unknown parameter and send their measurements to a fusion center. A fraction of the agents may be Byzantine agents that send arbitrary measurements to the fusion center, and the goal of the fusion center is to correctly infer the value of the unknown parameter even in the presence of Byzantine agents. In contrast to~\cite{Varshney1, Varshney2, FusionCenterEstimation1, FusionCenterDetection2}, this paper studies \textit{fully distributed} estimation: our setup does not use a fusion center.

Our previous work studies attacks against centralized cyber-physical systems -- we address the impact of information on detecting attacks~\cite{ChenTAC} and develop optimal attack strategies for adversaries~\cite{ChenAttack1, ChenAttack2}. Existing work in cyber-physical security also examines large-scale networked systems. Reference~\cite{MoGrid} provides examples of centralized cyber-physical attacks against power grids. Unlike~\cite{ChenTAC, ChenAttack1, ChenAttack2, MoGrid}, this paper studies resilience in a distributed setup instead of a centralized setup.

For fully distributed resilient algorithms, i.e., algorithms that do not depend on a fusion center, prior work has focused on average consensus~\cite{ResilientConsensus, LeBlanc1, Pasqualetti1},  fault detection in control systems~\cite{DistributedFD}, and general function computation~\cite{Sundaram1}. In~\cite{ResilientConsensus} and~\cite{LeBlanc1}, the authors consider the average consensus problem for scalar values in the presence of misbehaving agents. The authors provide iterative distributed consensus algorithms in which, during each iteration, agents ignore a subset of messages received from their neighbors. Reference~\cite{LeBlanc1} provides necessary and sufficient conditions, based on the number of adversarial agents and the topology of the inter-agent communication network, under which non-adversarial agents achieve consensus. The algorithms proposed in~\cite{ResilientConsensus} and~\cite{LeBlanc1} require each agent to have only local knowledge of the network topology, i.e., each agent needs to know its own neighbors in the network but does not need to know the entire network structure.  

Reference~\cite{Pasqualetti1} provides another method of dealing with adversarial agents in distributed consensus. The authors of~\cite{Pasqualetti1} construct fault detection and identification (FDI) filters for distributed consensus by leveraging knowledge of the network structure and algorithm dynamics (i.e., the process by which individual agents update their estimates and exchange information with neighbors). Unlike the methods of~\cite{ResilientConsensus} and~\cite{LeBlanc1}, the FDI filters proposed in~\cite{Pasqualetti1} depend on agents having knowledge of the entire network structure. Similarly, the authors of~\cite{DistributedFD} study distributed fault detection filters for interconnected dynamical systems (i.e., filters to detect whether or not a perturbing control signal has been applied to the system) that also depend on agents having knowledge of the entire network structure. 

When individual agents know the network structure, reference~\cite{Sundaram1} goes beyond average consensus and provides algorithms to calculate arbitrary functions of the agents' initial values in the presence of adversaries. In general, the resilience of the algorithms presented in~\cite{ResilientConsensus, LeBlanc1, Pasqualetti1, Sundaram1, DistributedFD} to adversarial agents depends on the topology of the inter-agent communication network. In this paper, we study distributed estimation instead of distributed consensus\footnote{We emphasize that, with distributed estimation, the agents make new measurements at every time instant, which contrasts with consensus where the data available is only the initial data, and no new measurements are available.} (like in~\cite{ResilientConsensus, LeBlanc1, Pasqualetti1, Sundaram1, DistributedFD}), and we propose an algorithm that does not require individual agents to know the entire network structure (unlike the algorithms presented in~\cite{Pasqualetti1, DistributedFD, Sundaram1}).

\subsection{Summary of Contributions}
In contrast with previous work on consensus, where no observations are involved, this paper studies resilient distributed estimation, where, at each time step, agents make their own observations, update their estimates, and exchange information with their neighbors. 
Reference~\cite{LeBlanc2} also considers attack resilient distributed parameter estimation, in which a group of agents attempt to estimate a \textit{scalar} parameter subject to some agents acting maliciously. In~\cite{LeBlanc2}, some agents have perfect information and know the value of the parameter before the estimation process, while other agents know a noisy version of the parameter. In this paper, unlike in~\cite{LeBlanc2}, no agent knows the value of the parameter before the estimation process.

We develop the Flag Raising Distributed Estimation ($\mathcal{FRDE}$) algorithm that simultaneously performs adversary detection and parameter estimation. The $\mathcal{FRDE}$ algorithm is a consensus+innovations estimator (see~\cite{Kar3}). In each iteration of $\mathcal{FRDE}$, each (normally-behaving) agent, to perform parameter estimation, updates its estimate based on its previous estimate, its sensor measurement (innovations), and the estimates of its neighbors (consensus), and, to perform adversary detection, checks the Euclidean distance between its own estimate and the estimates of its neighbors, reporting that an attack has occurred if the distance exceeds a certain (adaptive) threshold. 

Under the $\mathcal{FRDE}$ algorithm, if the network of normally-behaving agents is connected and the concatenation of their observations is globally observable\footnote{We formally define global observability in Section~\ref{sect: background}.}, then, either, the normally behaving agents will either asymptotically detect the presence of adversarial agents or their local estimates will be almost surely (a.s.) consistent (converge to true parameter value with probability one). That is, a set of adversarial agents that mask their presence in order not to alert others to their attack cannot simultaneously prevent the normally behaving agents' estimates from converging to the true value of the parameter. The key necessary and sufficient condition is global observability of the connected normally behaving agents. As long as this condition holds, consistency is almost surely guaranteed. Conversely, if the normally-behaving agents lack global observability, then, it is possible for adversarial agents to disrupt the estimation algorithm while simultaneously avoiding detection. 


The rest of this paper is organized as follows. In Section~\ref{sect: background}, we review background from spectral graph theory, present the sensing and communication model for the network of agents, and state the assumptions regarding adversarial agents. Section~\ref{sect: main} describes the $\mathcal{FRDE}$ algorithm, and Section~\ref{sect: analysis} analyzes its resilience to adversarial agents. We find an upper bound on the false alarm probability of the algorithm. The probability of missed detection depends on the behavior of the adversarial agents. When normally behaving agents lack global observability, adversarial agents may simultaneously evade detection and prevent local estimates from converging to the parameter. When normally behaving agents are connected and globally observable, then, they will either asymptotically detect the adversarial agents, or, in the case of no attack detection, their estimates will be a.s. consistent. We provide numerical examples of the $\mathcal{FRDE}$ algorithm in Section~\ref{sect: examples} and conclude in Section~\ref{sect: conclusion}. 

\section{Background}\label{sect: background}
\subsection{Notation}\label{sect: notation}
Let $\mathbb{R}^k$ denote the $k$ dimensional Euclidean space, $I_k$ the $k$ by $k$ identity matrix, $\mathbf{1}_k$ and $\mathbf{0}_k$ the column vectors of ones and zeros in $\mathbb{R}^k$, respectively. The operator $\left\lVert \cdot \right\rVert$ is the Euclidean $\ell_{2}$ norm when applied to vectors and  the corresponding induced norm when applied to matrices. For a matrix $M$, let $\mathscr{N}\left( M \right)$ be its null space and  $M^\dagger$ its Moore-Penrose pseudoinverse. The Kronecker product of $A$ and $B$ is $A \otimes B$. $M \succeq 0$ or $M \succ 0$ denotes that $M$ is positive semidefinite or positive definite, respectively.

Below, we consider a network of $N$ agents defined by a communication graph $G$. An undirected graph is noted as $G = \left(V, E\right)$, where $V = \left\{1, \dots, N\right\}$ is the set of vertices and  \[E = \left\{ (n, l) | \text{ $\exists$ an edge between vertex $n$ and vertex $l$} \right\}\] is the set of edges. We consider only simple graphs, i.e., no self loops nor multiple edges. The neighborhood of a vertex $n$ is
\[ \Omega_n = \left\{ l \in V \vert \left(n, l \right) \in E \right\}.\]
The degree of a vertex is $d_n = \left\vert \Omega_n \right\vert$. The degree matrix of $G$ is $D = \diag\left( d_1, \dots, d_N \right)$. The structure of G is described by the symmetric adjacency matrix $A = \left[ A_{nl} \right]$, where $A_{nl} = 1$ if $\left(n, l\right) \in E$ and $A_{nl} = 0$, otherwise. The positive semidefinite matrix $L = D-A$ is the graph Laplacian. The eigenvalues of $L$ can be ordered as $0 = \lambda_1 (L) \leq \dots \leq \lambda_N(L)$, and $\mathbf{1}_N$ is the eigenvector associated with $\lambda_1(L)$. For a connected graph $G$, $\lambda_2 (L) > 0$. For further details on spectral graph theory, see~\cite{Spectral, ModernGraph}.

For a subset of vertices $\mathcal{X} = \left\{ i_1, \dots, i_{\left\lvert \mathcal{X} \right \rvert} \right\} \subset V$, let $G_{\mathcal{X}} = \left( \mathcal{X}, E_{\mathcal{X}} \right)$, where
\[ E_{\mathcal{X}} = \left\{ (n, l) \in E | n \in \mathcal{X} \text{ and } l \in \mathcal{X} \right\}, \]
is the subgraph induced by $\mathcal{X}$. We say that a subset of vertices $\mathcal{X}$ is connected if $G_{\mathcal{X}}$ is connected. Let $L_{\mathcal{X}}$ denote the Laplacian of $G_{\mathcal{X}}$. For a vertex $n \in \mathcal{X}$ and for a subset of vertices $\mathcal{Y}$ disjoint from $\mathcal{X}$, let $\sigma_{n, \mathcal{Y}} = \left\vert \Omega_n \cap \mathcal{Y} \right\rvert$ denote the number of neighbors of $n$ in $\mathcal{Y}$. Let $\Sigma_{\mathcal{X}, \mathcal{Y}} = \diag \left( \sigma_{i_1, \mathcal{Y}}, \dots, \sigma_{i_{\left\lvert \mathcal{X} \right\rvert}, \mathcal{Y}} \right)$.

In this paper, we assume that all random objects are defined on a common probability space $\left(\Omega, \mathcal{F} \right)$. Let $\mathbb{P}\left(\cdot\right)$ and $\mathbb{E}\left[ \cdot \right]$ denote the probability and expectation operators, respectively. The abbreviation a.s. means ``almost surely,'' i.e., everywhere except on a set of measure $0$. In this paper, consistency of a estimator refers to strong consistency: a consistent estimator produces a sequence of estimates that converges a.s. to the parameter of interest.

\subsection{Sensing and Communication Model}\label{sect: model}
Consider a network of $N$ agents (or nodes) defined by a communication graph $G = \left(V, E\right)$. Let $\theta^* \in \mathbb{R}^M$ be a deterministic (distributed) unknown parameter that is to be estimated by the $N$ agents. Agent $n$ makes a measurement
\begin{equation}\label{eqn: sensingModel}
	y_n(t) = H_n \theta^* + w_n(t).
\end{equation}
For example, in power grid state estimation, $\theta^*$ represents the voltages and phase angles at all of the buses in the network. The $H_n$ matrices model local sensing. Again, in power grid state estimation, the $H_n$ model sensors measuring local voltages and phase angles at each of the buses~\cite{Xie2}. 

We make the following assumptions about the measurement noise term $w_n(t)$ and the measurement matrix $H_n$.
\begin{assumption}\label{ass: measurement}
	At each agent~$n$, the measurement noise $w_n(t)$ is independently and identically distributed (i.i.d.) over time $t$ with mean $\mathbb{E} \left[ w_n(t) \right] = 0$ and covariance $\mathbb{E} \left[ w_n(t) w_n(t)^T \right] = \Sigma_n$. Across agents, the measurement noise is independent, i.e., $w_j(t)$ and $w_k(s)$ are independently distributed for $j \neq k$ and for all $s, t$. 
\end{assumption}
\noindent For agent $n$, let $\psi_n = \trace\left( \Sigma_n \right)$. 
\begin{assumption}\label{ass: measureMatrix}
	At each agent $n$, the measurement matrix $H_n$ satisfies
	\begin{equation}\label{eqn: maxMeasurement}
		\lambda_{\max} \left( H_n^T H_n \right) \leq 1.
	\end{equation}
\end{assumption}
\noindent Assumption~\ref{ass: measureMatrix} is without loss of generality, since, if, for some $n$, $\lambda_{\max}\left(H_n^T H_n\right) > 1$, then appropriately scaling the measurement, i.e., $y_n(t)  = c_n \left( H_n \theta^* + w_n(t) \right)$, we can get $\lambda_{\max}\left(c^2_n H_n^T H_n \right) \leq 1$. 

The goal of the network of agents is to recursively estimate $\theta^*$ from the individual observations $y_1(t), \dots, y_N(t)$, $t = 0, 1, \dots$, subject to the assumptions stated below.
\begin{assumption}\label{ass: connectivity}
	The graph $G$ is connected.
\end{assumption}
\noindent If $G$ is not connected, then we can separately consider each connected component of $G$.
\begin{specification}\label{ass: localKnowledge}
Agent $n$ knows only its own local sensing model (i.e., each agent knows only its own $H_n$), its neighborhood $\Omega_n$, and its local observation $y_n(t)$. 
\end{specification}

\begin{specification}\label{ass: localExchange}
	Agent $n$ exchanges information only with other agents in $\Omega_n$.
\end{specification}

\begin{definition}\label{def: globalObservable}
	Let $\mathcal{X} = \left\{ n_1, \dots, n_{\left\vert \mathcal{X} \right\vert} \right\}$ be a subset of (connected) agents ($\mathcal{X} \subseteq V$) and $G_{\mathcal{X}}$ its induced subgraph. We say that $G_{\mathcal{X}}$ is globally observable if the matrix $\sum_{n \in \mathcal{X}} H_n^TH_n$ is invertible.
\end{definition}
\noindent Global observability is a necessary condition for a centralized estimator to be consistent. So, it is natural to assume it for distributed estimation (e.g.,~\cite{Kar2}) as we state next.
\begin{assumption}\label{ass: globalObservable}
	The graph of all agents, $G$, is globally observable.
\end{assumption}
\noindent The individual agents are not assumed to be locally observable, i.e., $H_n$ may have low column rank.

\begin{assumption}\label{ass: compactness}
	The parameter $\theta^*$ belongs to a non-empty compact set $\Theta$, where
\begin{equation}\label{eqn: thetaSetDef}
	\Theta = \left\{ \theta \in \mathbb{R}^M \vert \left\lVert \theta \right\rVert \leq \eta \right\}.
\end{equation}
	Each agent $n$ knows the value of $\eta$. 
\end{assumption}

\noindent In practical settings for distributed parameter estimation, such as distributed state estimation in the power grid~\cite{Xie2} and temperature estimation over wireless sensor networks~\cite{Kar2}, the parameter to be estimated does not take an arbitrary value but rather takes a value from a closed, bounded set determined by physical laws.

\subsection{Threat Model}
The agents use an iterative distributed message-passing protocol, to be specified shortly, to estimate the parameter $\theta^*$. Some agents in the network are adversarial and attempt to disrupt the estimation procedure. We partition the set of all agents $V$ into a set of adversarial agents $\mathcal{A}$ and a set of normal agents $\mathcal{N} = V\setminus \mathcal{A}$. Agent $n$ is adversarial, $n \in \mathcal{A}$, if, for some $t = 0, 1, \dots$, it deviates from the distributed protocol (to be introduced soon). Our goal in this paper is to develop secure protocols for distributed parameter estimation.
\begin{specification}\label{ass: powerfulAttacker}
	The adversarial agents know the structure of $G$, the true value of the parameter $\theta^*$, the members of the sets $\mathcal{A}$ and $\mathcal{N}$, may communicate with all other adversarial agents to launch powerful attacks, and send arbitrary messages to their neighbors.
\end{specification}
\noindent  Normally behaving agents, $n \in \mathcal{N}$, do not initially know whether other agents are normally behaving or adversarial. That is, the agents $n\in\mathcal{N}$ do not know the members of $\mathcal{N}$ and $\mathcal{A}$.

Our attack model differs from the model presented in~\cite{LMSCluster1, LMSCluster2}. In these references, the adversarial agents act as though they observe a different parameter than the true parameter $\theta^*$ and share this distorted information with their neighbors. The intruders are not required to know the true parameter or the size of the network. In our  paper, the adversarial agents know the true parameter and the size of the network and may send different messages to each of its neighbors.  It is in this sense, that ours is a worst-case scenario of adversarial agents.  Our attack model includes the attack described in~\cite{LMSCluster1, LMSCluster2}, i.e. in this paper, adversarial agents may attack the network by acting as though they observe a different parameter than $\theta^*$, but they are not required to attack in this particular way.


For the purpose of analysis, only,  we make the following assumptions regarding the sets of normally-behaving and adversarial agents.
\begin{assumption}\label{ass: normalConnected}
	The induced subgraph $\mathcal{G}_{\mathcal{N}}$ of the normally behaving agents is connected.
\end{assumption}
\noindent In practice, Assumption~\ref{ass: normalConnected} may not be true; adversarial agents may split the normally behaving agents into several connected components. In such a situation, our analysis applies to each connected component of normally behaving agents. For the remainder of this paper, unless otherwise stated, we make Assumption~\ref{ass: normalConnected}.

\section{Distributed Estimation and Adversary Detection}\label{sect: main}
In this section, we present the Flag Raising Distributed Estimation ($\mathcal{FRDE}$) algorithm, under which agents simultaneously perform parameter estimation and adversary detection.

\subsection{$\mathcal{FRDE}$ Algorithm}\label{sect: algorithm} Each iteration of $\mathcal{FRDE}$ follows three main steps: \begin{enumerate*} \item Message Passing, \item Estimate Update, and \item Adversary Detection \end{enumerate*}.
 At each $t = 0, 1, \dots$, each agent $n$ computes an estimate $x_n(t)$ of the parameter $\theta^*$  and a flag value $\pi_n(t)$. The flag $\pi_n(t)$ takes values of either ``Attack'' or ``No Attack,'' depending if agent $n$ detects an attack. We say that the distributed algorithm has detected an adversarial agent if for some $t$ and $n$, $\pi_n(t) = \text{Attack}$, and we say that the distributed algorithm has not detected an adversarial agent if for all $t$ and all $n$, $\pi_n(t) = \text{No Attack}$.\footnote{By Assumption~\ref{ass: normalConnected}, $G_{\mathcal{N}}$ is connected, so, once a single agent $n \in \mathcal{N}$ detects an attack, it communicates this detection to all other agents in $\mathcal{N}$.} The flag $\pi_n (t) = \text{``Attack''}$ indicates that there exists an adversarial agent in the network but does not identify which agent(s) is (are) adversarial. 

We initialize the $\mathcal{FRDE}$ algorithm as follows. At time $t=0$, each agent $n \in \mathcal{N}$ sets its estimate as
\begin{equation}\label{eqn: agentInitialEstimate}
	x_n(0) = 0,
\end{equation}
and its initial flag value as
\begin{equation}\label{eqn: agentInitialFlag}
	\pi_n(0)= \text{No Attack}.
\end{equation}
Note that, as a result of~\eqref{eqn: agentInitialEstimate}, all of the initial estimates of normally behaving agents $n$ satisfy $x_n(0) \in \Theta$. 

\underline{Message Passing}: For all $t = 0, 1, 2, \dots$, the normally behaving agents $n\in\mathcal{N}$ follow the message generation rule
\begin{equation}\label{eqn: normalMessage}
	m_{n, l}^t = x_n(t),
\end{equation}
and send $m_{n, l}^t$ to each of its neighbors.
 
\underline{Estimate Update}: Each agent maintains a running average of its local measurement 
\begin{equation}\label{eqn: normalAverageMeasurement}
	\begin{split}
		\overline{y}_n(t) &= {t \over t+1}\overline{y}_n(t-1) + {1 \over t+1}y_n(t), \\
		\overline{y}_n(0) &= y_n(0).
	\end{split}
\end{equation}
Note that
\begin{equation}\label{eqn: averageMeasurementEquiv}
	\overline{y}_n(t) = H_n \theta^* + {1 \over t+1}\sum_{j=0}^t w_n(j).
\end{equation}
Agent $n \in \mathcal{N}$ follows the consensus+innovations estimation update rule
\begin{equation}\label{eqn: normalUpdate}
\begin{split}
	x_n(t+1) =& x_n(t) - \beta \sum_{l \in \Omega_n} \left( x_n(t) - m_{l, n}^t \right) \\ &+ \alpha H_n^T \left(\overline{y}_n(t) - H_n x_n(t) \right),
\end{split}
\end{equation}
where $\alpha$ and $\beta$ are positive constants to be specified shortly.

\underline{Adversary Detection}: Agent $n \in \mathcal{N}$ updates its flag by the rule
\begin{equation}\label{eqn: normalFlagUpdate}
	\pi_n(t+1) = \!\!\left\{\begin{array}{ll}\text{Attack}, & \pi_n(t) = \text{Attack}, \text{ or } \\
		& \exists l, \left \lVert x_n(t) - m_{l, n}^t \right\rVert > \gamma_t \\
		\text{No Attack}, & \text{Otherwise} \end{array}\right.,
\end{equation}
where $\gamma_t$ is a time-varying parameter to be specified shortly. We assume that no adversarial agent purposefully reports an attack (i.e.,  for all $n \in \mathcal{A}$, and for all $t$, $\pi_n(t) = \text{No Attack}$), since the adversarial agents want to avoid being detected. The threshold parameter $\gamma_t$ follows the recursion
\begin{equation}\label{eqn: gammaDef}
\begin{split}
	\gamma_{t+1} &= \left(1 - r_1 \right)\gamma_t + \alpha\frac{2K}{(t+1)^{\tau}}, \\
	\gamma_0 &= 2\eta \sqrt{N},
\end{split}
\end{equation}
where $K > 0$, $0 < \tau < \frac{1}{2}$, and $0 < r_1 \leq 1$ are parameters to be specified shortly. We describe how to select parameters $\alpha$, $\beta$, $r_1$, $K$, and $\tau$ in Section~\ref{sect: param}. We require $\alpha$, $\beta$, and $r_1$ to satisfy the following conditions:
\begin{equation}\label{eqn: maxLambda}
	\lambda_{\max} \left( J_{\beta, \alpha}\right) \leq 1,
\end{equation}
\begin{equation}\label{eqn: rCond}
	0 <	r_1 \leq \lambda_{\min} \left( J_{\beta, \alpha} \right),
\end{equation}
where $J_{\beta, \alpha} = \beta \left(L \otimes I_M \right) + \alpha D^T_H D_H$, $D_H = \blkdiag\left(H_1, \dots, H_N \right)$, and we recall $M$ is the dimension of~$\theta^*$.

Intuitively, the definition of $\gamma_t$ in~\eqref{eqn: gammaDef} and conditions~\eqref{eqn: maxLambda} and~\eqref{eqn: rCond} relate to the performance of $\mathcal{FRDE}$ as follows. Following~\eqref{eqn: gammaDef}, the threshold $\gamma_t$ decays over time. If we choose parameters to satisfy~\eqref{eqn: maxLambda} and~\eqref{eqn: rCond}, then, in the absence of adversaries, the local estimation error decays in a manner similar to the threshold $\gamma_t$. Moreover, in the absence of adversaries, the Euclidean distance between any agent's local estimate and any of its received messages is upper bounded by a constant factor times the local estimation error. Thus, if this upper bound is ever violated, i.e., if the Euclidean distance between the local estimate and any received message exceeds the threshold, the agent reports the presence of an adversary. We provide details on how conditions~\eqref{eqn: maxLambda} and~\eqref{eqn: rCond} relate to the performance of $\mathcal{FRDE}$ in Section~\ref{sect: analysis}.


The estimate update rule in~\eqref{eqn: normalUpdate} is exactly the same as the distributed parameter estimation algorithm provided in~\cite{Xie2}, and, if there are no adversarial nodes, the distributed algorithm described in this paper behaves exactly as the distributed algorithm provided in~\cite{Xie2}.\footnote{In~\cite{Xie2}, the step sizes $\alpha$ and $\beta$ decay over time to cope with measurement noise. This paper considers constant $\alpha$ and $\beta$.} Here, we emphasize that our main contributions are the distributed attack detection update described by equations~\eqref{eqn: normalFlagUpdate} and~\eqref{eqn: gammaDef} and the analysis of $\mathcal{FRDE}$ under adversarial activities. Unlike the existing literature~\cite{Kar1, Xie2, Kar2} on consensus+innovation parameter estimation, which assumes all agents behave normally, our attack detection method allows the distributed parameter estimation method to operate even in the presence of adversarial agents. For DLMS algorithms, references~\cite{LMSCluster1} and~\cite{LMSCluster2} study the effect of adversarial agents who behave as though they observe a different parameter than the true parameter.

\subsection{Parameter Selection}\label{sect: param}
This subsection describes how to select parameters $\alpha, \beta, r_1, K$ and $\tau$. for the $\mathcal{FRDE}$ algorithm. Parameters $K$ and $\tau$ may take any values that satisfy $K > 0$ and $0 < \tau < {1 \over 2}$. Section~\ref{sect: analysis} describes how the choices of $K$ and $\tau$ affect the performance of the algorithm. We describe two procedures to select the parameters $\alpha$, $\beta$, and $r_1$. We assume that, during the setup phase of $\mathcal{FRDE}$ (for parameter selection), all agents behave normally. 

Procedure 1 requires centralized knowledge of the entire network structure and all sensing matrices $H_1, \dots, H_N$ during the setup phase of $\mathcal{FRDE}$. Once the parameters $\alpha$, $\beta$, and $r_1$ are chosen, besides these three parameters, individual agents only need local information (its own measurement and the messages of its neighbors) to execute the $\mathcal{FRDE}$ algorithm. They do not need to know the structure of the entire network and the  sensing matrices of the other agents. Procedure 2 describes how to select parameters knowing only that the network $G$ is connected (Assumption~\ref{ass: connectivity}) and globally observable (Assumption~\ref{ass: globalObservable}). That is, Procedure 2 does not require exact knowledge of the network structure and sensing matrices.

We now establish that Procedures 1 and 2 satisfy~\eqref{eqn: maxLambda} and~\eqref{eqn: rCond}. We start with the following Lemma
\begin{lemma}\label{lem: JPD}
	For any $\alpha, \beta > 0$, the matrix $J_{\beta, \alpha} = \beta \left( L \otimes I_M \right) + \alpha D_H^T D_H$ is positive definite.
\end{lemma}
\noindent The proof of Lemma~\ref{lem: JPD} is found in the appendix. The first procedure is as follows.

\noindent\textbf{Procedure 1:}
\begin{enumerate}
	\item Choose auxiliary parameters $\widehat{\alpha}$ and $\widehat{\beta}$ to be positive.
	\item Set \begin{equation}\label{eqn: 1ab} \alpha = \frac{\widehat{\alpha}}{\lambda_{\max} \left( J_{\widehat{\beta}, \widehat{\alpha}}\right)}, \quad \beta= \frac{\widehat{\beta}}{\lambda_{\max} \left( J_{\widehat{\beta}, \widehat{\alpha}}\right)}.\end{equation}
	\item Set $r_1 = \lambda_{\min} \left( J_{\beta, \alpha}\right).$
\end{enumerate}

\begin{lemma}\label{lem: proc1Conditions}
	The parameters $\alpha, \beta,$ and $r_1$ selected using Procedure 1 lead to $\lambda_{\max} \left( J_{\beta, \alpha} \right) \leq 1$ and $0 < r_1 \leq \lambda_{\min}\left( J_{\beta, \alpha} \right)$. 
\end{lemma}
\noindent The proof of Lemma~\ref{lem: proc1Conditions} is found in the appendix.
Computing $\lambda_{\max} \left( J_{\widehat{\beta}, \widehat{\alpha}} \right)$ in~\eqref{eqn: 1ab} requires knowledge of the network structure (the graph Laplacian, $L$) and all sensing matrices (the matrix $D_H^T D_H$). 

Knowledge of the network structure and all sensing matrices is a design cost associated with Procedure 1. Individual agents, however, do not need to know or store the network structure and sensing matrices; they only need to store the scalar parameters $\alpha, \beta,$ and $r_1$ to perform $\mathcal{FRDE}$. We can compute these parameters separately in the cloud and broadcast them to the agents, so that no agent needs to know the graph $G$ or the $H_n$ of other agents. Moreover, in practice, the network structure and sensing matrices are sparse. For example, in a wireless temperature sensor network, where $\theta^*$ represents a temperature field, every sensor measures a single component of $\theta^*$, and the resulting $H_n$ matrices are sparse (each $H_n$ matrix has exactly one nonzero component). The sparsity of the network and $H_n$ matrices means that this information can be efficiently shared amongst all of the agents, which, in practice, reduces the design cost associated with Procedure~1.

We present a second procedure that only requires knowledge of $\lambda_{\max} (L)$, the second smallest eigenvalue, $\lambda_2 (L)$, of the graph Laplacian $L$, and the minimum eigenvalue of the matrix $\mathcal{G} = \frac{1}{N} \sum_{n = 1}^N H_n^T H_n$. 

\noindent\textbf{Procedure 2:}
\begin{enumerate}
	\item Choose $\kappa_1 > \frac{1}{\lambda_2 (L)} \left(\lambda_{\min} (\mathcal{G}) + 2 \sqrt{4 - \lambda_{\min} (\mathcal{G})}\right).$
	\item Set $\alpha = \left(\kappa_1 \lambda_{\max} \left(L\right) + N \right)^{-1}$.
	\item Set $\beta = \alpha \kappa_1$. 
	\item Set $r_1 = \alpha \left(\lambda_{\min} (\mathcal{G}) - \frac{4\lambda_{\min} (\mathcal{G})}{\sqrt{4 \lambda_{\min} (\mathcal{G}) - \left( \lambda_{\min}(\mathcal{G}) - \lambda_2(L)\kappa_1 \right)^2}} \right).$
\end{enumerate}

\noindent If we do not know $\lambda_{\max} (L)$ and $\lambda_{2} (L)$, we can use bounds for these quantities that depend only on the number of agents $N$~\cite{Graphs}. In particular, in steps 1) and 4), we can replace $\lambda_{2} (L)$ with the lower bound $\lambda_2 (L) \geq \frac{4}{N^2}$, and, in step 2), we can replace $\lambda_{\max} (L)$ with the upper bound $\lambda_{\max}(L) \leq {N}$. 

\begin{lemma}\label{lem: proc2Conditions}
	The parameters $\alpha, \beta,$ and $r_1$ selected using Procedure 2 lead to $\lambda_{\max} \left( J_{\beta, \alpha} \right) \leq 1$ and $0 < r_1 \leq \lambda_{\min}\left( J_{\beta, \alpha} \right)$. 
\end{lemma}
\noindent The proof of Lemma~\ref{lem: proc2Conditions} is found in the appendix. 

\section{Performance Analysis of $\mathcal{FRDE}$}\label{sect: analysis}
This section studies the performance of the $\mathcal{FRDE}$ algorithm and its resilience to adversarial agents. Under the $\mathcal{FRDE}$ algorithm, normally behaving agents will either achieve a.s. consistency in their local estimates or they will asymptotically report the presence of adversaries. The necessary and sufficient condition for this performance is the global observability of the (connected) network of normally behaving agents. In the case that all agents behave normally (i.e., no adversaries), we show that the probability of any agent incorrectly reporting an adversary can be made arbitrarily small. That is, we find an upper bound for the algorithm's probability of false alarm. 


Let $\mathbf{x}_t = \left[\begin{array}{ccc} x_1(t)^T & \cdots & x_N(t)^T \end{array} \right]^T,$ and define the global estimation error
\begin{equation}\label{eqn: rConProof3}
	\mathbf{e}_t = \mathbf{x}_t - \left( 1_N \otimes \theta^* \right).
\end{equation}
We study the behavior of $\mathbf{e}_t$ and of the flag variable $\pi_n(t)$ over time. In this section, we assume the parameters $\alpha, \beta,$ and $r_1$ are chosen to satisfy~\eqref{eqn: maxLambda} and~\eqref{eqn: rCond}. 


\subsection{Performance with No Adversarial Agents}
We find an upper bound for the false alarm probability of $\mathcal{FRDE}$ (i.e., the probability that $\mathcal{FRDE}$ declares that there is an adversarial agent under the scenario that all agents behave normally) and analyze the behavior of local estimates $x_n(t)$ when there are no adversarial agents.

\begin{theorem}\label{thm: regularConvergence}
	If there are no adversarial agents ($\mathcal{A} = \emptyset$), then, under the $\mathcal{FRDE}$ algorithm, for all $n \in V$, we have
\begin{equation}\label{eqn: regularConvergence}
	\mathbb{P} \left( \lim_{t \rightarrow \infty} \left(t+1\right)^{\tau_0} \left\lVert x_n(t) - \theta^* \right\rVert = 0 \right) = 1,
\end{equation}
for every $0 \leq \tau_0 < {1 \over 2}$.
Moreover, the false alarm probability, $P_{FA}$, satisfies
\begin{equation}\label{eqn: falseAlarm}
	P_{FA} = \mathbb{P} \left(\exists n \in V, t \geq 0 :  \pi_n(t) =  \text{Attack} \right) \leq \frac{\Psi \zeta(\tau)}{K^2},
\end{equation}
where $\Psi = \sum_{n = 1}^N \trace\left(\Sigma_n\right)$ and $\zeta(\tau) = \sum_{j = 1}^\infty \frac{1}{j^{2\left(1-\tau\right)}}$.
\end{theorem}
\noindent Theorem~\ref{thm: regularConvergence} states that, in the absence of adversaries, local estimates are strongly consistent, converging almost surely to $\theta^*$. Equation~\eqref{eqn: falseAlarm} bounds the probability that \textit{any} agent $n$ at \textit{any} time $t \geq 0$ raises an alarm flag. Therefore,~\eqref{eqn: falseAlarm} bounds the false alarm probability of $\mathcal{FRDE}$. 

To decrease the false alarm probability, one should choose larger values of the parameter $K$ and smaller values of the parameter $\tau$. For any $0 < \tau < {1 \over 2}$, we can make the false alarm probability arbitrarily small by choosing large enough $K$. The upper bound provided by~\eqref{eqn: falseAlarm} is conservative. To achieve a low false alarm rate in practice, we may not need as large a value of $K$ as dictated by~\eqref{eqn: falseAlarm}. We illustrate the difference between the upper bound on false alarm probability and the algorithm's empirical false alarm rate through numerical examples in Section~\ref{sect: examples}. Choosing larger $K$ and smaller $\tau$ results in the threshold $\gamma_t$ decaying more slowly (i.e., larger $\gamma_t$). Having a larger $\gamma_t$ allows the adversarial agents to send more malicious messages (messages that deviate more from true parameter) while evading detection. We illustrate the trade off between the false alarm probability and the magnitude of the threshold $\gamma_t$ through numerical examples in Section~\ref{sect: examples}. 

The proof of Theorem~\ref{thm: regularConvergence} requires several intermediate results. We will use the following lemma to determine the effect of measurement noise on the estimation process.

\begin{lemma} [Lemma 5 in~\cite{Kar1}]\label{lem: tvSystem}Consider the scalar, time varying system
\begin{equation}\label{eqn: tvSystem}
	v_{t+1} = \left(1 - p_1(t) \right) v_t + p_2(t),
\end{equation}
where
\begin{equation}\label{eqn: tvRates}
	p_1(t) = \frac{c_1}{(t+1)^{\delta_1}}, \: p_2(t) = \frac{c_2}{(t+1)^{\delta_2}},
\end{equation}
$c_1, c_2 > 0$, $0 \leq \delta_1 \leq 1$, and $\delta_2 \geq 0$. If $\delta_1 < 1$ and $\delta_2 > \delta_1$, we have
\begin{equation}\label{eqn: tvLimit}
	\lim_{t \rightarrow \infty} (t+1)^{\delta_0} v_t  = 0,
\end{equation}
for every $0 \leq \delta_0 < \delta_2 - \delta_1$. If $\delta_1 = 1$ and $\delta_2 > \delta_1$,~\eqref{eqn: tvLimit} holds if, in addition, $\delta_0 < c_1$. 
\end{lemma}
\noindent The proof of Lemma~\ref{lem: tvSystem} may be found in~\cite{Kar1} and~\cite{Kar3}. 

The following result characterizes the behavior of time-averaged measurement noise.
\begin{lemma}\label{lem: timeAveragedNoise}
	Let $w_0, w_1, w_2, \dots$ be i.i.d. (vector) random variables with mean $\mathbb{E} \left[ w_t \right] = 0$ and finite covariance $\mathbb{E} \left[ w_t w_t^T\right] = \Sigma$. Define $m_t$ as the mean of $w_0, \dots, w_t$, i.e.,
\begin{equation}\label{eqn: averageNoise}
	m_t = \frac{1}{t+1} \sum_{j = 0}^t w_j.
\end{equation}
Then, we have
\begin{equation}\label{eqn: averageNoiseConvergence}
	\mathbb{P} \left(\lim_{t\rightarrow \infty} (t+1)^{\tau_0} \left\lVert m_t \right\rVert = 0 \right) = 1,\end{equation}
\begin{equation}\label{eqn: averageNoiseSup}
	\mathbb{P}\left(\sup_{t \geq 0} \left\lVert m_t \right\rVert >  \frac{K}{(t+1)^{\tau_0}}\right) \leq \frac{\trace\left(\Sigma\right)\zeta(\tau_0)}{K^2},
\end{equation}
for $0 \leq \tau_0 < \frac{1}{2}$, where $\zeta(\tau_0) = \sum_{j = 1}^\infty \frac{1}{j^{2\left(1-\tau_0\right)}}$.
\end{lemma}
\noindent The proof of Lemma~\ref{lem: timeAveragedNoise} may be found in the appendix. 

We now prove Theorem~\ref{thm: regularConvergence}.
\begin{IEEEproof}[Proof (Theorem~\ref{thm: regularConvergence})]
First, we study the a.s. convergence of the estimates and prove~\eqref{eqn: regularConvergence}. Let \[\overline{\mathbf{y}}_t = \left[\begin{array}{ccc} \overline{y}_1(t)^T & \cdots & \overline{y}_N(t)^T \end{array} \right]^T,\] and let \[\overline{\mathbf{w}}_t = \left[\begin{array}{ccc} \overline{w}_1(t)^T & \cdots & \overline{w}_N(t)^T \end{array} \right]^T,\] where $\overline{w}_n(t) = \frac{1}{t+1} \sum_{j=0}^t w_j$. From~\eqref{eqn: normalUpdate}, we have
\begin{equation}\label{eqn: rConProof1}
	\mathbf{x}_{t+1} = \mathbf{x}_t -\beta \left( L \otimes I_M \right) \mathbf{x}_t + \alpha D_H^T \left( \overline{\mathbf{y}}_t - D_H \mathbf{x}_t \right).
\end{equation}
Performing algebraic manipulation and noting, from~\eqref{eqn: sensingModel}, that \[D_H^T \overline{\mathbf{y}}_t = D_H^T D_H \left( 1_N \otimes \theta^* \right) + D_H^T\overline{\mathbf{w}}_t,\] we have
\begin{equation}\label{eqn: rConProof2}
	\mathbf{x}_{t+1} = \left(I_{NM} - J_{\beta, \alpha} \right) \mathbf{x}_t + J_{\beta, \alpha}\left(1_N \otimes \theta^* \right) + \alpha D_H^T \overline{\mathbf{w}}_t,
\end{equation}
where~\eqref{eqn: rConProof2} follows from~\eqref{eqn: rConProof1} because, from the properties of $L$, we have $\left(L \otimes I_M \right)\left(1_N \otimes \theta^* \right) = 0$. 

From~\eqref{eqn: rConProof2}, we have that the dynamics of $\mathbf{e}_t$ are
\begin{equation}\label{eqn: rConProof4}
	\mathbf{e}_{t+1} = \left(I_{NM} - J_{\beta, \alpha} \right) \mathbf{e}_{t} + \alpha D_H^T \overline{\mathbf{w}}_t. 
\end{equation}
Conditions~\eqref{eqn: maxLambda} and~\eqref{eqn: rCond} state that $\lambda_{\max} \left(J_{\beta, \alpha} \right) \leq 1$ and $\lambda_{\min} \left( J_{\beta, \alpha} \right) > 0$. Thus the matrix $I_{NM} - J_{\beta, \alpha}$ is positive semi-definite with
\begin{equation}\label{eqn: rConProof5}
	\lambda_{\max} \left( I_{NM} - J_{\beta, \alpha} \right) = 1 - \lambda_{\min} \left(J_{\beta, \alpha} \right) < 1. 
\end{equation}
By Assumption~\ref{ass: measureMatrix}, we have $\lambda_{\max} \left( H_n^T H_n \right) \leq 1$, which means that $\left\lVert D_H^T \right \rVert \leq 1$. Since $\left\lVert I_{NM} - J_{\beta, \alpha} \right\rVert = \lambda_{\max} \left(I_{NM} - J_{\beta, \alpha} \right) < 1$, we have, from~\eqref{eqn: rConProof4}, that
\begin{align}\label{eqn: rConProof6}
	\left\lVert \mathbf{e}_{t+1} \right \rVert &\leq \left \lVert I_{NM} - J_{\beta, \alpha} \right \rVert \left \lVert \mathbf{e}_t \right \rVert+ \alpha \left \lVert D_H^T \overline{\mathbf{w}}_t \right \rVert \\
	& \leq \left( 1 - \lambda_{\min} \left( J_{\beta, \alpha} \right) \right) \left \lVert \mathbf{e}_t \right \rVert + \alpha \left\lVert \overline{\mathbf{w}}_t \right \rVert. \label{eqn: rConProof7}
\end{align}
By construction, $\overline{\mathbf{w}}_t$ falls under the purview of Lemma~\ref{lem: timeAveragedNoise}, which means that
\begin{align}\label{eqn: rConProof8}
	\mathbb{P} \left( \lim_{t \rightarrow \infty} (t+1)^{\tau_1} \left\lVert \overline{\mathbf{w}}_t \right\rVert = 0 \right) &= 1,
\end{align}
for every $0 \leq \tau_1 < \frac{1}{2}$. As a consequence of~\eqref{eqn: rConProof8}, there exists finite $K_{\overline{\mathbf{w}}}$ such that
\begin{equation}\label{eqn: rConProof9}
	\mathbb{P} \left ( \left\lVert \overline{\mathbf{w}}_t \right\rVert \leq {K_{\overline{\mathbf{w}}} \over (t+1)^{\tau_1}} \right) = 1,
\end{equation}
for every $0 \leq \tau_1 < \frac{1}{2}$. Thus, almost surely, we have, from~\eqref{eqn: rConProof7},
\begin{equation}\label{eqn: rConProof10}
	\left\lVert \mathbf{e}_{t+1} \right \rVert \leq \left( 1 - \lambda_{\min} \left( J_{\beta, \alpha} \right) \right) \left \lVert \mathbf{e}_t \right \rVert + \frac{\alpha K_{\overline{\mathbf{w}}}}{(t+1)^{\tau_1}}.
\end{equation}
The relationship above falls under the purview of Lemma~\ref{lem: tvSystem}, which means that we have
\begin{equation}\label{eqn: rConProof11}
	\mathbb{P} \left(\lim_{t\rightarrow \infty} (t+1)^{\tau_0} \left \lVert \mathbf{e}_t \right\rVert = 0 \right) = 1,
\end{equation}
for all $0 \leq \tau_0 < \tau_1$. Since $\left\lVert x_n(t) - \theta^* \right \rVert \leq \left \lVert \mathbf{e}_t \right \rVert$, by taking $\tau_1$ arbitrarily close to $1 \over 2$, equation~\eqref{eqn: rConProof11} establishes~\eqref{eqn: regularConvergence}. 

Second, we bound the false alarm probability of $\mathcal{FRDE}$. The $\mathcal{FRDE}$ algorithm raises a false alarm if, in the absence of adversaries, for any $t \geq 0$, any $n \in V$, and any $l \in \Omega_n$, we have $\left\lVert x_n(t) - x_l(t) \right \rVert > \gamma_t$. By the triangle inequality, we have
\begin{align}\label{eqn: falseAlarms1}
	\left\lVert x_n(t) - x_l(t) \right \rVert &\leq \left\lVert x_n(t) - \theta^* \right\rVert + \left \lVert x_l(t) - \theta^* \right\rVert, \\&\leq 2 \left\lVert \mathbf{e}_t \right \rVert.\label{eqn: falseAlarms2}
\end{align}
Thus, if there is a false alarm, then, necessarily, for some $t \geq 0$, $\left\lVert \mathbf{e}_t \right \rVert > {\gamma_t \over 2}$. We find an upper bound on the probability that $\left\lVert \mathbf{e}_t \right \rVert > {\gamma_t \over 2}$. 

As an intermediate step, recall, from above, that $\overline{\mathbf{w}}_t$ falls under the purview of Lemma~\ref{lem: timeAveragedNoise}. Thus, for any $K > 0$, $0 \leq \tau < \frac{1}{2}$, we have
\begin{equation}\label{eqn: falseAlarms1a}
	\mathbb{P} \left(\sup_{t \geq 0} \left\lVert \overline{\mathbf{w}}_t \right \rVert >  {K \over (t+1)^{\tau}} \right) \leq \frac{\Psi \zeta(\tau)}{K^2},
\end{equation}
where $\Psi = \sum_{n = 1}^N \trace\left(\Sigma_n\right)$ and $\zeta(\tau) = \sum_{j = 1}^\infty \frac{1}{j^{2\left(1-\tau\right)}}$. Consider the set of sample paths \begin{equation}\label{eqn: samplePaths1} \left\{ \omega \in \Omega : \sup_{t \geq 0} \left\lVert \overline{\mathbf{w}}_t \right \rVert \leq {K \over (t+1)^{\tau}} \right\}.\end{equation}
By~\eqref{eqn: falseAlarms1a}, such a set has probability greater than $1 - \frac{\Psi \zeta(\tau)}{K^2}$, and we use induction to show that, on this set, $\left\lVert \mathbf{e}_t \right \rVert \leq {\gamma_t \over 2}$ for all $t \geq 0$.

In the base case, for $t = 0$, since $x_n(0) = 0$ and $\left\lVert \theta^* \right \rVert \leq \eta$, we have $\left\lVert \mathbf{e}_t \right \rVert \leq \eta \sqrt{N} = {\gamma_t \over 2}$. In the induction step, we assume that $\left\lVert \mathbf{e}_t \right \rVert \leq {\gamma_t \over 2}$ and show that $\left \lVert \mathbf{e}_{t+1} \right \rVert \leq {\gamma_{t+1} \over 2}$. By~\eqref{eqn: rCond}, we have $\lambda_{\min} \left(J_{\beta, \alpha} \right) \geq r_1$. Since we only consider sample paths on the set defined by~\eqref{eqn: samplePaths1}, we have $\left\lVert \overline{\mathbf{w}}_t \right \rVert \leq {K \over (t+1)^{\tau}}$. Then, from~\eqref{eqn: rConProof7}, we have 
\begin{align}
	\left\lVert \mathbf{e}_t \right \rVert &\leq \left(1 - r_1 \right)\left\lVert \mathbf{e} \right \rVert + \alpha \frac{K}{(t+1)^{\tau}}, \label{eqn: falseAlarms2} \\
	&\leq \left(1-r_1 \right) {\gamma_t \over 2} + \alpha \frac{K}{(t+1)^{\tau}} = {\gamma_{t+1} \over 2}.\label{eqn: falseAlarms3}
\end{align}
From~\eqref{eqn: falseAlarms3}, we conclude that, on the set defined by~\eqref{eqn: samplePaths1}, i.e., with probability greater than $1 - \frac{\Psi \zeta(\tau)}{K^2}$, $\left\lVert \mathbf{e}_t \right \rVert \leq {\gamma_t \over 2}$ for all $t \geq 0$. Then, the false alarm rate of $\mathcal{FRDE}$ is upper bounded by $\frac{\Psi \zeta(\tau)}{K^2}$, which establishes~\eqref{eqn: falseAlarm}. 
\end{IEEEproof}

\subsection{Performance with Adversarial Agents}
When there are adversarial agents present, the performance of the $\mathcal{FRDE}$ algorithm depends on the strength of the adversarial and normally behaving agents. Specifically, the algorithm's performance depends on the global observability of the normally behaving agents. If the network of  normally behaving agents is not globally observable, then, the adversarial agents may disrupt the estimation process while evading detection. We say that a set of adversarial agents $\mathcal{A}$ evades detection by $\mathcal{FRDE}$ if the probability of detecting $\mathcal{A}$ is no greater than the false alarm probability of $\mathcal{FRDE}$. 

\begin{proposition}\label{prop: undetectableAttack}
	Let the parameter $\theta^*$ satisfy $\left\lVert \theta^* \right\rVert < \eta$. If the network $G_{\mathcal{N}}$ of the normally behaving agents $\mathcal{N}$ is not globally observable, then it is possible for the set of adversarial agents $\mathcal{A}$ to perform an attack (i.e., to send messages to their neighbors) such that, for some $\overline{\theta} \neq \theta^*$, all $n \in \mathcal{N}$, and every $0 \leq \tau_0 < {1 \over 2}$,
\begin{equation}\label{eqn: fakeConvergence}
	\mathbb{P} \left( \lim_{t \rightarrow \infty} \left(t+1\right)^{\tau_0} \left\lVert x_n(t) - \overline{\theta} \right\rVert = 0 \right) = 1,
\end{equation}
and
\begin{equation}\label{eqn: fakeFalseAlarm}
	\mathbb{P} \left(\exists n \in \mathcal{N}, t \geq 0 :  \pi_n(t) =  \text{Attack} \right) \leq P_{FA}.
\end{equation}
where $P_{FA}$, given in~\eqref{eqn: falseAlarm}, is the false alarm rate of $\mathcal{FRDE}$.
\end{proposition}
\noindent Proposition~\ref{prop: undetectableAttack} states that, if the normally behaving agents are not globally observable, then the adversarial agents may attack the algorithm in a way that simultaneously evades detection and causes all local estimates to converge to a wrong parameter, $\overline{\theta}$, almost surely. Proposition~\ref{prop: undetectableAttack} holds even if the set of adversarial agents, $\mathcal{A}$, does not induce a connected subgraph~$G_{\mathcal{A}}$. 
\begin{IEEEproof}
	To prove Proposition~\ref{prop: undetectableAttack}, we design an attack that simultaneously evades the distributed attack detection algorithm and causes the normally behaving agents' estimates to almost surely converge to an incorrect value. Let $\mathcal{N} = \left\{ n_1, \dots, n_{| \mathcal{N} | } \right\}$, and let $\mathcal{H}_{\mathcal{N}} = \left[\begin{array}{ccc} H_{n_1}^T & \cdots & H_{n_{|\mathcal{N}|}} \end{array} \right]^T$. Since the network of normally behaving agents is not globally observable, we have
\[
	\mathcal{H}_{\mathcal{N}}^T \mathcal{H}_{\mathcal{N}} = \sum_{n\in \mathcal{N}} H_n^T H_n
\]
is not invertible, which means that there exists a nonzero $\overline{\mu} \in \mathbb{R}^M$ such that
\begin{equation}\label{eqn: undetectableAttack1}
	\mathcal{H}_{\mathcal{N}} \overline{\mu} = 0.
\end{equation}
Choose $\overline{\mu} \in \mathscr{N} \left( \mathcal{H}_{\mathcal{N}} \right)$ such that
\begin{equation}\label{eqn: undetectableAttack2}
	\left\lVert \theta^* + \overline{\mu} \right\rVert \leq \eta. 
\end{equation}
Such a choice of $\overline{\mu}$ always exists, since, $\mathscr{N} \left( \mathcal{H}_{\mathcal{N}} \right)$ is a subspace, $\left\lVert \theta^* \right \rVert < \eta$ (from the statement of Proposition~\ref{prop: undetectableAttack}), and, from the triangle inequality, we have $\left\lVert \theta^* + \overline{\mu} \right\rVert \leq \left \lVert \theta^* \right\rVert + \left\lVert \overline{\mu} \right\rVert$. 

Let $\overline{\theta} = \theta^*  + \overline{\mu}$. Recall, from Specification~\ref{ass: powerfulAttacker}, that all adversarial agents know the value of the parameter $\theta^*$. Let all adversarial agents $n \in \mathcal{A}$ participate in the distributed estimation algorithm as though the true parameter is $\overline{\theta}$. That is, an agent $n \in \mathcal{A}$, behaves as though its sensor measurement~is
\begin{equation}\label{eqn: undetectableAttack3}
	\overline{y}_n = H_n \overline{\theta}, 
\end{equation}
and, otherwise, it follows the initial estimate generation, message generation, and estimate update rules of the distributed estimation algorithm, given by equations~\eqref{eqn: agentInitialEstimate},~\eqref{eqn: normalMessage}, and~\eqref{eqn: normalUpdate}, respectively. Note that this is the same attack strategy as described in~\cite{LMSCluster1, LMSCluster2}.

Consider the scenario in which the true parameter to be estimated is $\overline{\theta}$. For all $n \in \mathcal{N}$, we have $H_n \overline{\theta} = H_n \theta^*$, since, by definition, $\overline{\theta} - \theta^* \in \mathscr{N} \left( \mathcal{H}_{\mathcal{N}} \right)$. That is, in the scenario that the true parameter to be estimated is $\overline{\theta}$, the agents in $\mathcal{N}$ make the same sensor measurements (up to measurement noise) as in the scenario that the true parameter is $\theta^*$. Thus, the case in which the true parameter is $\theta^*$ but the set of adversarial agents behaves as though the true parameter is $\overline{\theta}$ is equivalent to the case in which the true parameter is $\overline{\theta}$ and all agents behave normally. Thus, the adversaries probability of being detected is no greater than the false alarm rate of $\mathcal{FRDE}$. Equation~\eqref{eqn: fakeConvergence} follows as a consequence of Theorem~\ref{thm: regularConvergence}.
\end{IEEEproof}

We now consider the case when the set of (connected) uncompromised agents is globally observable. One of two events must occur: either \begin{enumerate*} \item there exists some uncompromised agent $n \in \mathcal{N}$ that raises an alarm flag ($\pi_n(t) = \text{Attack}$), or \item no uncompromised agent ever raises an alarm flag (i.e., for all $n \in N$ and for all $t \geq 0$, $\pi_n(t) = \text{No Attack}$).\end{enumerate*} If event 1) occurs, then, the $\mathcal{FRDE}$ successfully detects the presence of an adversarial agent. If event 2) occurs, then, the $\mathcal{FRDE}$ algorithm has a missed detection. The following theorem states that, in the case that $\mathcal{FRDE}$ misses a detection (event 2), the local estimates of normally behaving agents are consistent.
 
\begin{theorem}\label{thm: resilientOperation}
	Let the set of normally behaving agents $\mathcal{N}$ be connected, and let $G_{\mathcal{N}}$ be globally observable. If, for all $n \in \mathcal{N}$ and for all $t = 0, 1, \dots$, we have $\pi_n(t) = \text{No Attack}$, then, under $\mathcal{FRDE}$, for all $n \in \mathcal{N}$, we have
\begin{equation}\label{eqn: resilientConvergence}
	\mathbb{P} \left( \lim_{t \rightarrow \infty} \left(t+1\right)^{\tau_0} \left\lVert x_n(t) - \theta^* \right\rVert = 0 \right) = 1,
\end{equation}
for every $0 \leq \tau_0 < {\tau}$.
\end{theorem}
\noindent Theorem~\ref{thm: resilientOperation} states that when the normally behaving agents are connected and their models are globally observable, if the adversarial agents are undetected (for all times $t \geq 0$), then, almost surely, all normally behaving agents' local estimates converge asymptotically to $\theta^*$. 

For the remainder of this subsection, without loss of generality, let $\mathcal{N} = \left\{1, \dots, \left\vert\mathcal{N} \right\vert \right\}$, and let $\mathcal{A} = \left\{ \left\vert \mathcal{N} \right\vert+1, \dots, N \right\}$.\footnote{Although the normally behaving agents $n \in \mathcal{N}$ are not initially aware of the members of $\mathcal{N}$ and $\mathcal{A}$, for purposes of analysis, we can relabel the agents arbitrarily without loss of generality.} Let $\mathbf{x}_t^{\mathcal{N}} = \left[\begin{array}{ccc} x_1(t)^T & \cdots & x_{\left\vert \mathcal{N} \right \vert} (t)^T \end{array} \right]^T,$ and define the estimation error of the normally behaving agents as 
\begin{equation}
	\mathbf{e}_t^{\mathcal{N}} = \mathbf{x}_{t}^{\mathcal{N}} - \left(1_{\left\vert \mathcal{N} \right \vert} \otimes \theta^* \right).
\end{equation}
To prove Theorem~\ref{thm: resilientOperation}, we study the behavior of $\mathbf{e}_t^{\mathcal{N}}$ over time, and we require the following Lemma, the proof of which is found in the appendix. 

\begin{lemma}\label{lem: LPartition}
Let $G = (V, E)$ be a graph, and, for $n \in V$, let $H_n$ be the sensing matrix associated with agent $n$. Let $\mathcal{X} = \left\{ i_1, \dots, i_{\left\vert \mathcal{X} \right\vert} \right\} \subset V$ be a globally observable subset of agent models ( i.e., the matrix $\sum_{i \in \mathcal{X}} H_i^T H_i$ is invertible) that induces a connected subgraph $G_{\mathcal{X}}$ of $G$. Suppose that $\alpha, \beta > 0$ are chosen such that \[J_{\beta, \alpha} = \beta \left(L \otimes I_M \right) + \alpha D_H^T D_H\] is positive definite and satisfies $\lambda_{\max} \left(J_{\beta, \alpha} \right) \leq 1$. Then, the matrix
\begin{equation}\label{eqn: JPartition}
	J^{\mathcal{X}}_{\beta, \alpha} = \beta\left( L_{\mathcal{X}} \otimes I_M \right) + \alpha {D_H^{\mathcal{X}}}^T {D_H^{\mathcal{X}}},
\end{equation}
where $L_{\mathcal{X}}$ is the graph Laplacian of $G_{\mathcal{X}}$ and $D_H^{\mathcal{X}} = \blkdiag\left( H_{i_1}, \dots, H_{i_{\left\vert \mathcal{X} \right\vert}}\right)$ 
is also positive definite and satisfies $\lambda_{\max} \left( J^{\mathcal{X}}_{\beta, \alpha} \right) \leq 1$. 
\end{lemma}

We now prove Theorem~\ref{thm: resilientOperation}.
\begin{IEEEproof}[Proof (Theorem~\ref{thm: resilientOperation})]
Consider the update equation of a normally behaving $n \in \mathcal{N}$.  For any node $n \in V$, we can partition the neighborhood of $n$ into $\Omega_n^{\mathcal{N}} = \Omega_n \cap \mathcal{N}$ and $\Omega_n^{\mathcal{A}} = \Omega_n \cap \mathcal{A}$. We have
\begin{equation}\label{eqn: resilient1}
	\begin{split}
		&x_n(t+1) = x_n(t) - \beta \sum_{l \in \Omega_n^{\mathcal{N}}} \left( x_n(t) - x_l(t) \right) - \\
		&\quad \beta \sum_{l \in \Omega_n^{\mathcal{A}}} \left(x_n(t) - m_{l, n}^t \right) + \alpha H_n^T \left(\overline{y}_n(t) - H_n x_n(t) \right).
	\end{split}
\end{equation}
Equation~\eqref{eqn: resilient1} follows from~\eqref{eqn: normalMessage} and~\eqref{eqn: normalUpdate} because all normally behaving agents follow the prescribed message generation rule while all adversarial agents may send arbitrary messages. The condition that, for all $n \in \mathcal{N}$ and for all $t$, we have $\pi_n (t) = \text{No Attack}$, implies that for all $l \in \Omega_n$, we have \[\left\lVert x_n(t) - m_{l, n}^t \right \rVert \leq \gamma_t.\]
Thus, we can rewrite~\eqref{eqn: resilient1} as
\begin{equation}\label{eqn: resilient2}
		\begin{split}
		&x_n(t+1) = x_n(t) - \beta \sum_{l \in \Omega_n^{\mathcal{N}}} \left( x_n(t) - x_l(t) \right) + \\
		&\quad \alpha H_n^T \left(\overline{y}_n(t) - H_n x_n(t) \right) + \delta_n(t),
	\end{split}
\end{equation}
where $\delta_n(t)$ is a bounded disturbance vector that satisfies 
\begin{equation}\label{eqn: resilient3}
	\left\lVert \delta_n(t) \right \rVert \leq \beta \left\lvert \mathcal{A} \right\rvert \gamma_t.
\end{equation}


Let \[ \Delta_t^{\mathcal{N}} = \left[\begin{array}{ccc} \delta_1(t)^T & \cdots & \delta_{\left \vert \mathcal{N} \right \vert}(t)^T \end{array} \right]^T.\] From~\eqref{eqn: resilient3}, we have
\begin{equation}\label{eqn: resilient5}
	\left\lVert \Delta_t^{\mathcal{N}} \right\rVert \leq \beta \left\vert \mathcal{A} \right \vert\sqrt{\left\lvert \mathcal{N} \right\rvert}\gamma_t.
\end{equation}
Using~\eqref{eqn: resilient2}, we compute the dynamics of $\mathbf{x}_t^{\mathcal{N}}$ as
\begin{equation}\label{eqn: resilient6}
\begin{split}
	\mathbf{x}_{t+1}^{\mathcal{N}} =& \left(I_{\left \lvert \mathcal{N} \right\rvert M}- \beta \left( L_{\mathcal{N}} \otimes I_M \right) \right) \mathbf{x}_t^{\mathcal{N}} + \Delta^{\mathcal{N}}_t \\
	&+ {D_H^{\mathcal{N}}}^T \left( \overline{\mathbf{y}}_t^{\mathcal{N}} - \alpha {D_H^{\mathcal{N}}} \mathbf{x}_t^{\mathcal{N}} \right),
\end{split}
\end{equation}
where $L_{\mathcal{N}}$ is the Laplacian of the induced subgraph $\mathcal{G}_{\mathcal{N}}$, $D_H^{\mathcal{N}} = \blkdiag \left( H_1, \dots, H_{\left\vert \mathcal{N}\right\vert} \right)$, and $\overline{\mathbf{y}}_t^{\mathcal{N}} = \left[\begin{array}{ccc} \overline{y}_1(t)^T & \cdots & \overline{y}_{\left\vert \mathcal{N} \right\vert}(t)^T \end{array} \right]^T$. 

Following the same algebraic manipulations as in the proof of Theorem~\ref{thm: regularConvergence}, we have
\begin{equation}\label{eqn: resilient7}
	\mathbf{e}_{t+1}^{\mathcal{N}} = \left(I_{\left\vert\mathcal{N}\right\vert M } - J_{\beta, \alpha}^{\mathcal{N}} \right) \mathbf{e}_{t}^{\mathcal{N}} + \Delta_{t}^{\mathcal{N}} + \alpha {D_H^{\mathcal{N}}}^T \overline{\mathbf{w}}_t^{\mathcal{N}}, 
\end{equation}
where $J_{\beta, \alpha}^{\mathcal{N}} = \beta \left(L_{\mathcal{N}} \otimes I_M \right) + \alpha \left( {D_H^{\mathcal{N}}}^T \mathcal{D}_H^{\mathcal{N}}\right)$ and $\overline{\mathbf{w}}_t^{\mathcal{N}} = \left[\begin{array}{ccc} \overline{w}_1(t)^T & \cdots & \overline{w}_{\left\vert \mathcal{N} \right\vert}(t)^T \end{array} \right]^T$. By the triangle inequality we have
\begin{align}
\begin{split}\label{eqn: resilientConv1}
	\left \lVert \mathbf{e}_{t+1}^{\mathcal{N}} \right\rVert &\leq \left \lVert I_{\left\vert\mathcal{N}\right\vert M } - J_{\beta, \alpha}^{\mathcal{N}} \right \rVert \left \lVert \mathbf{e}_{t}^{\mathcal{N}}  \right \rVert + \alpha \left \lVert {D_H^{\mathcal{N}}}^T \overline{\mathbf{w}}_t^{\mathcal{N}} \right \rVert +  \\ 
	& \quad\beta \left\vert \mathcal{A} \right\vert \sqrt{\left\lvert \mathcal{N} \right \rvert} \gamma_t,
\end{split}
\end{align}
Note that $\overline{\mathbf{w}}_t^{\mathcal{N}}$ falls under the purview of Lemma~\ref{lem: timeAveragedNoise}. Thus, we have
\begin{equation}\label{eqn: normalNoise1}
	\mathbb{P} \left( \lim_{t\rightarrow\infty} (t+1)^{\tau_1} \left\lVert \overline{\mathbf{w}}_t^{\mathcal{N}} \right \rVert = 0 \right) = 1,
\end{equation}
for every $0 \leq \tau_1 < {1 \over 2}$, and, as a consequence of~\eqref{eqn: normalNoise1}, there exists finite $K^{\mathcal{N}}_{\overline{\mathbf{w}}}$ such that
\begin{equation}\label{eqn: normalNoise2}
	\mathbb{P} \left( \left\lVert \overline{\mathbf{w}}_t^{\mathcal{N}}\right\rVert \leq \frac{K^{\mathcal{N}}_{\overline{\mathbf{w}}}}{(t+1)^{\tau_1}} \right) = 1.
\end{equation}
Also, note that $\gamma_t$, as defined in~\eqref{eqn: gammaDef}, falls under the purview of Lemma~\ref{lem: tvSystem}, so we have
\begin{equation}\label{eqn: gammaEvolve}
	\lim_{t\rightarrow\infty} (t+1)^{\tau_2} \gamma_t = 0,
\end{equation}
for every $0 \leq \tau_2 < \tau$. As a consequence of~\eqref{eqn: gammaEvolve}, there exists finite $\Gamma > 0$ such that
\begin{equation}\label{eqn: gammaBound}
	\gamma_t \leq \frac{\Gamma}{(t+1)^{\tau_2}}.
\end{equation}

Substituting~\eqref{eqn: normalNoise2} and~\eqref{eqn: gammaBound} into~\eqref{eqn: resilientConv1}, we have, almost surely, that
\begin{align}\label{eqn: resilientConv2}
	\begin{split}
	\left\lVert \mathbf{e}_{t+1}^{\mathcal{N}} \right \rVert&\leq \left( 1 - \lambda_{\min} \left(J_{\beta, \alpha}^{\mathcal{N}} \right)\right)\left\lVert \mathbf{e}_t^{\mathcal{N}} \right \rVert + \alpha \frac{K_{\overline{\mathbf{w}}}^{\mathcal{N}}}{(t+1)^{\tau_1}} + \\
	& \quad \beta \frac{\left\vert \mathcal{A} \right\vert \sqrt{\left\lvert \mathcal{N} \right \rvert}\Gamma}{(t+1)^{\tau_2}}.
\end{split}
\end{align}
Taking $\tau_1 \geq \tau_2$ and $\tau_2$ arbitrarily close to $\tau$ yields
\begin{align}\label{eqn: resilientConv3}
	\left\lVert \mathbf{e}_{t+1}^{\mathcal{N}} \right \rVert \leq \left( 1 - \lambda_{\min} \left(J_{\beta, \alpha}^{\mathcal{N}} \right)\right)\left\lVert \mathbf{e}_t^{\mathcal{N}} \right \rVert + \frac{c_3}{(t+1)^{\tau}},
\end{align}
for some $c_3 > 0$. The recurrence relation in~\eqref{eqn: resilientConv3} falls under the purview of Lemma~\ref{lem: tvSystem}, which means that
\begin{equation}\label{eqn: resilientConv4}
	\lim_{t \rightarrow \infty} (t+1)^{\tau_0} \left\lVert \mathbf{e}_t^{\mathcal{N}} \right \rVert = 0.
\end{equation}
For any $n \in \mathcal{N}$, $\left\lVert x_n(t) - \theta^* \right \rVert \leq \left\lVert \mathbf{e}_t^{\mathcal{N}} \right \rVert$, so~\eqref{eqn: resilientConv4} establishes~\eqref{eqn: resilientConvergence}.
\end{IEEEproof}

\section{Numerical Examples}\label{sect: examples}

We demonstrate the performance of the $\mathcal{FRDE}$ algorithm. The motivation for the numerical examples is as follows. Consider a network of mobile agents, for example, robots, whose goal is to determine and arrive at an (initially) unknown target location. The agents are equipped with sensors to measure the target location, and all of them can communicate over a fixed communication graph $G$. Individual robots do not know the entire structure of $G$ and instead know only their local neighborhood in $G$. The normally behaving robots use the $\mathcal{FRDE}$ algorithm to estimate the location of the target from their collective measurements and report the presence of adversarial robots. The adversarial robots attempt to cause an error in the distributed estimation process while avoiding detection. 

\begin{figure}[h!]
	\centering
	\includegraphics[keepaspectratio = true, scale = .6]{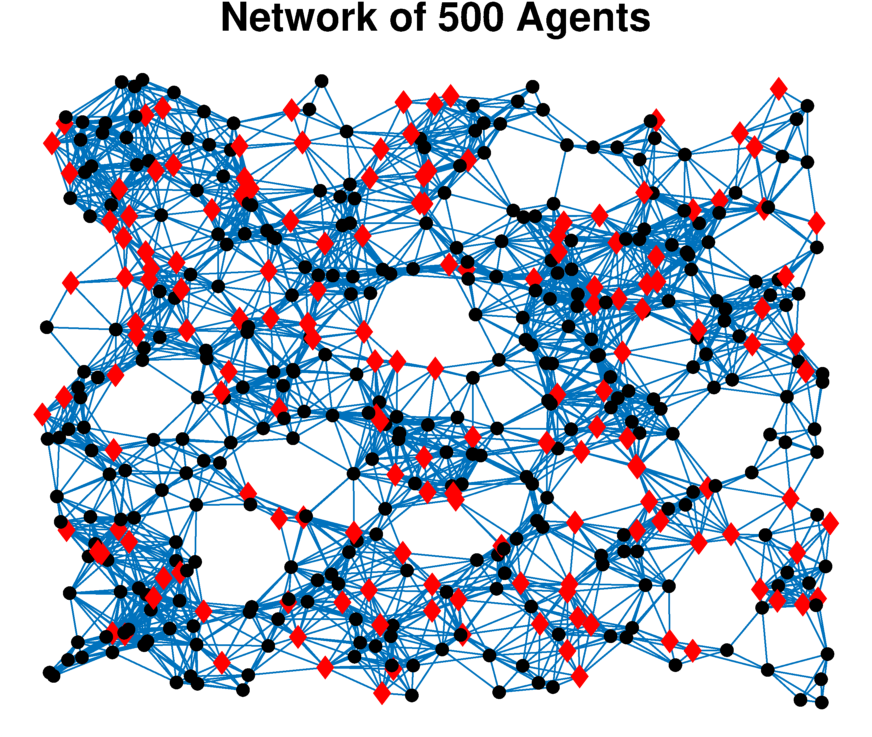}
	\caption{Communication Network of the 500 agents. Agents marked by black dots measure the $x$ and $y$ coordinates of the target. Agents marked by red diamonds measure the $z$ coordinate of the target.}
	\label{fig: network}
\end{figure}

We consider a network of $N = 500$ agents attempting to estimate the parameter $\theta^* \in \mathbb{R}^3$, which corresponds to the $x, y,$ and $z$ coordinates of the target. For all examples, we consider the same $\theta^*$: we choose $\theta^*$ (uniformly) at random from a sphere of radius $\eta = 500$ meters. The agents communicate over a random geometric network, given by Figure~\ref{fig: network}.
We place agents uniformly at random over a two dimensional square (100 meters by 100 meters) and place an edge between agents whose Euclidean distance is below 10 meters.

We randomly select $160$ agents, indicated by the red diamonds in Figure~\ref{fig: network}, to be equipped with sensors that measure the $z$ component of the target location. For such agents, we have 
\[H^{\text{Diamond}}_n = \left[\begin{array}{ccc} 0 & 0 & 1 \end{array} \right].\]The remaining agents, indicated by the black dots, are equipped with sensors that measure the $x$ and $y$ components of the target location. For such agents, we have \[H^{\text{Circle}}_n = \left[\begin{array}{ccc} 1 & 0 & 0 \\ 0 & 1 & 0 \end{array} \right]. \]Additive measurement noise affects all agents' sensors. The measurement noise for every agent $n$, $w_n(t)$, is an i.i.d. sequence of Gaussian random variables with mean $0$ and covariance $\sigma^2 I_{p_n}$, where $p_n$ is the dimension of the measurement $y_n(t)$. For our numerical examples, we use the covariance value $\sigma^2 = 60$. The local signal-to-noise ratio (SNR) is $11$~dB. 

The agents perform the algorithm $\mathcal{FRDE}$ to estimate $\theta^*$. For the numerical examples, we assign the value 0 to the $\text{No Attack}$ flag and the value $1$ to the $\text{Attack}$ flag. We use Procedure 1 to compute the following parameters for $\mathcal{FRDE}$: $\alpha = 3.1 \times 10^{-2}$, $\beta = 3.1 \times 10^{-2}$, $r_1 = 9.0 \times 10^{-3}$. These parameter choices satisfy conditions~\eqref{eqn: maxLambda} and~\eqref{eqn: rCond}. We also choose parameters $K = 4$ and $\tau = 0.40$.

We consider four different configurations of adversarial agents:
\begin{enumerate}
	\item \textbf{No adversarial agents:} All agents behave normally,
	\item \textbf{Strong adversarial agents:} All red diamond agents are adversarial. The remaining normally behaving agent models are globally unobservable. 
	\item \textbf{Disruptive weak adversarial agents:} Half of the red diamond agents are adversarial and perform a disruptive attack. That is, the adversarial agents attempt to compromise the consistency of the remaining agents' estimates. The remaining normally behaving agent models are globally observable.
	\item \textbf{Undisruptive weak adversarial agents:} Half of the red diamond agents are adversarial and perform a stealthy attack. That is, the adversarial agents attack the remaining agents so that no agent raises an alarm flag. The remaining normally behaving agent models are globally observable.
\end{enumerate}
For each configuration of the adversarial agents, we show how the local estimates and flag values evolve in time for a single execution of $\mathcal{FRDE}$.

\subsection{Local Estimate and Flag Value Evolution}
\begin{figure}[h!]
	\centering
	\includegraphics[keepaspectratio = true, scale = .6]{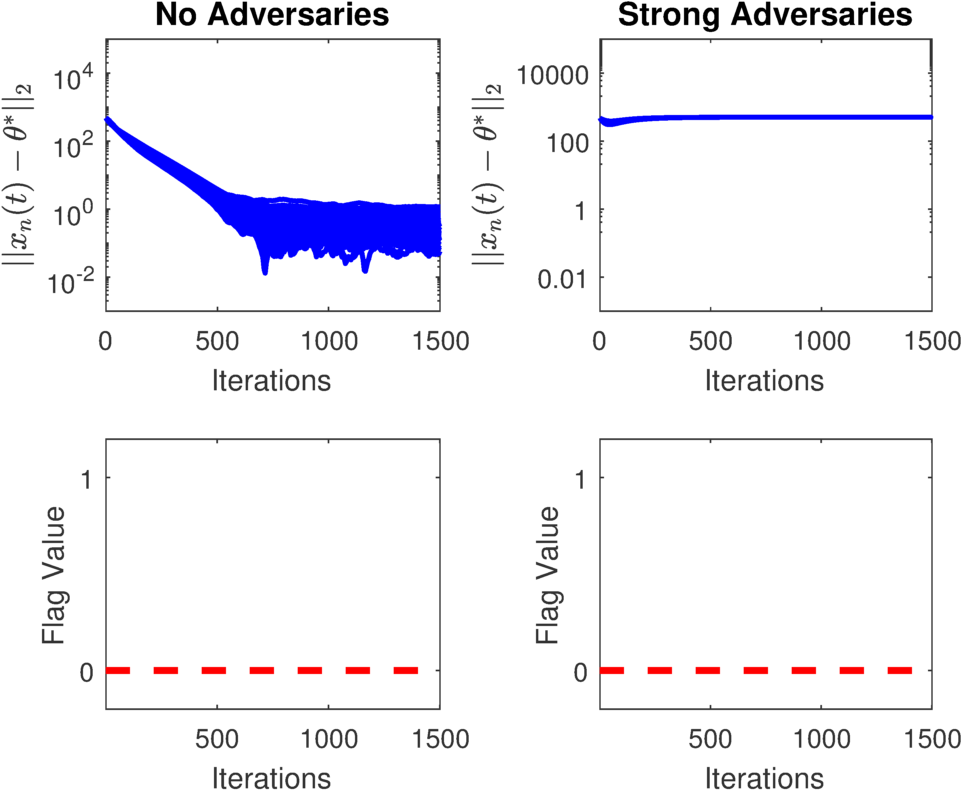}
	\caption{Performance of $\mathcal{FRDE}$ when there are no adversarial agents (left) and when all $160$ red diamond agents are adversarial (right). Top: Agent Estimation Errors. Bottom: Agent Flag Values.}
	\label{fig: noAttack}
\end{figure}

Figure~\ref{fig: noAttack} describes the performance of $\mathcal{FRDE}$ when all agents behave normally and when all red diamond agents are adversarial.
In the absence of adversarial agents, following Theorem~\ref{thm: regularConvergence}, the local estimates converge to $\theta^*$.

When all $160$ red diamond agents are adversarial, the normally behaving agents induce a connected subgraph but are not globally observable. Then, following Proposition~\ref{prop: undetectableAttack}, it is possible for the adversarial agents to simultaneously avoid detection and cause the normally behaving agents to estimate $\theta^*$ incorrectly. Specifically, the adversarial agents can behave as though the true parameter was $\theta^* + \overline{\mu}$, where $\overline{\mu}$ is an offset in the parameter's $z$ coordinate that satisfies $\left\lVert \theta^* + \mu \right \rVert \leq \eta$. Figure~\ref{fig: noAttack} shows that, following the attack strategy of Proposition~\ref{prop: undetectableAttack}, the adversarial agents can prevent the network of agents from converging to the correct estimate while remaining undetected.
\begin{figure}[h!]
	\centering
	\includegraphics[keepaspectratio = true, scale = .6]{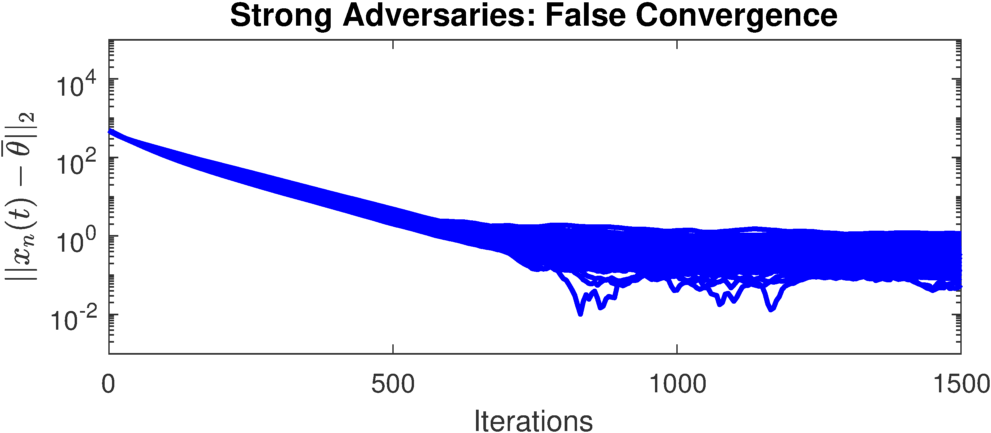}
	\caption{Performance of $\mathcal{FRDE}$ when all $160$ red diamond agents are adversarial: Convergence of estimates to $\theta^* + \overline{\mu}$.}
	\label{fig: falseConvergence}
\end{figure}
As shown by Figure~\ref{fig: falseConvergence}, under the attack of Proposition~\ref{prop: undetectableAttack}, the estimates of all agents converge to the incorrect parameter $\theta^* + \overline{\mu}$.

In the third and fourth numerical example, we randomly select $80$ out of the $160$ red diamond agents to be adversarial. In these examples, the network of normally behaving agents is connected and globally observable. First, we consider the case in which the $80$ adversarial agents behave as though the true parameter was $\theta^* + \overline{\mu}$. 
\begin{figure}[h!]
	\centering
	\includegraphics[keepaspectratio = true, scale = .6]{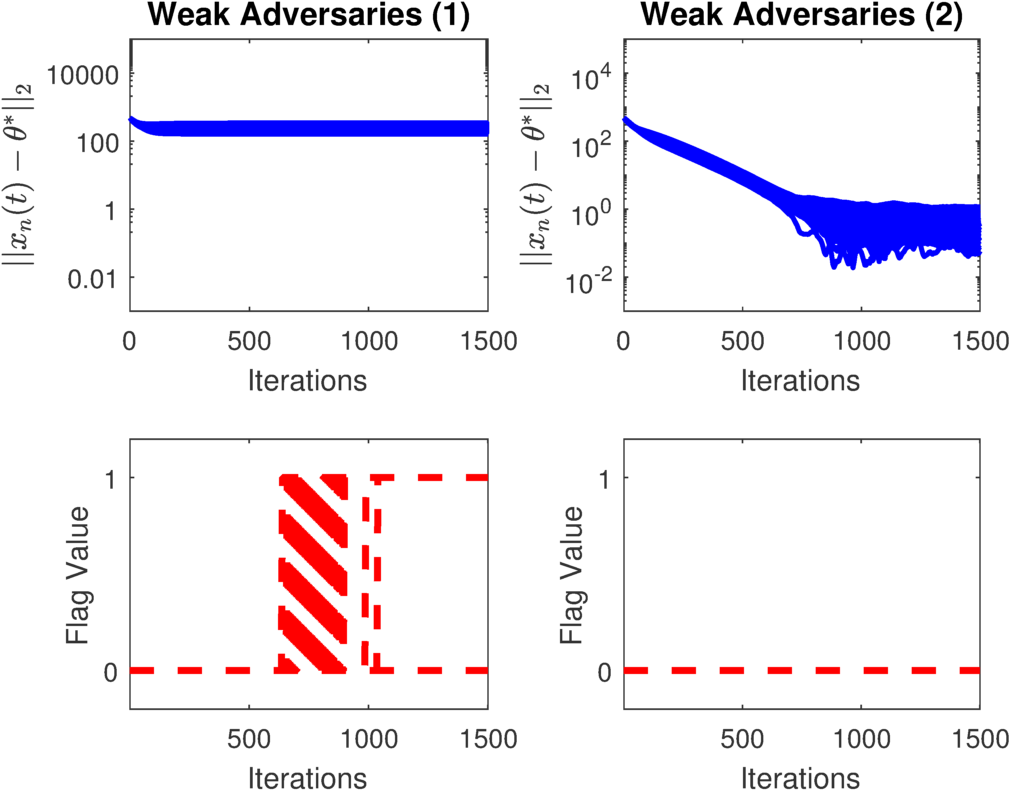}
	\caption{Performance of $\mathcal{FRDE}$ when $80$ out of $160$ red diamond agents are adversarial. Left: Disruptive Attack Right: Undisruptive attack. Top: Agent Estimation Errors. Bottom: Agent Flag Values.}
	\label{fig: unobservableAttack}
\end{figure}
Figure~\ref{fig: unobservableAttack} shows that the adversarial agents are able to prevent the normally behaving agents' estimates from converging to the correct value. The attack, however, causes the normally behaving agents to raise flags and indicate the presence of an adversarial agent. Thus, although the agents do not converge to the correct estimate, the $\mathcal{FRDE}$ algorithm allows the normally behaving agents to correctly detect the presence of an adversarial agent.

In the fourth numerical example, we again consider the case in which $80$ of the $160$ red diamond agents are adversarial. Although, as previously noted, the normally behaving agent models in this case are globally observable, it may still be possible for the adversarial agents to perform an undetectable attack. To avoid detection,  for any agent $n \in \mathcal{N}$, for any $l \in \Omega_n$, and for all iterations $t$, the estimate $x_n(t)$, must satisfy $\left\lVert x_n(t) - m_{l, n}^t \right\rVert \leq \gamma_t$. That is, to avoid detection, adversarial agents must attack the network in such a way that no agent's estimate deviates too far from the estimates of its neighbors.

Figure~\ref{fig: unobservableAttack} shows the effect of an undetectable attack when the normally behaving agent models are globally observable. Adversarial agents avoid detection, as no agent raises a flag indicating the presence of an adversarial agent, but the adversarial agents cannot prevent the normally behaving agents from converging to the correct estimate of $\theta^*$. The results of the third and fourth numerical examples verify Theorem~\ref{thm: resilientOperation}. If the network of normally behaving agents is connected globally observable, and, if the adversarial agents behave in an undetectable manner, the normally behaving agents' estimates converge to the parameter $\theta^*$.

\subsection{Performance Trade-offs}
In this subsection, we evaluate the performance trade-offs for $\mathcal{FRDE}$ between the false alarm probability and the deviation tolerated in adversarial messages. Specifically, for the network described by~Figure~\ref{fig: network}, we show how the false alarm probability bound in~\eqref{eqn: falseAlarm} and the evolution of the threshold $\gamma_t$ depend on the parameters $K$ and $\tau$ (for fixed $\alpha$, $\beta$, and $r_1$). It is desirable to have both small false aparm probability and small $\gamma_t$. Having a small $\gamma_t$ means that, in order to avoid detection, adversarial agents may not send messages that deviate too far from the receiving agents' estimate. Conversely, having a large $\gamma_t$ allows adversarial agents to send more malicious messages while evading detection.

\begin{figure}[h!]
	\centering
	\includegraphics[keepaspectratio = true, scale = .6]{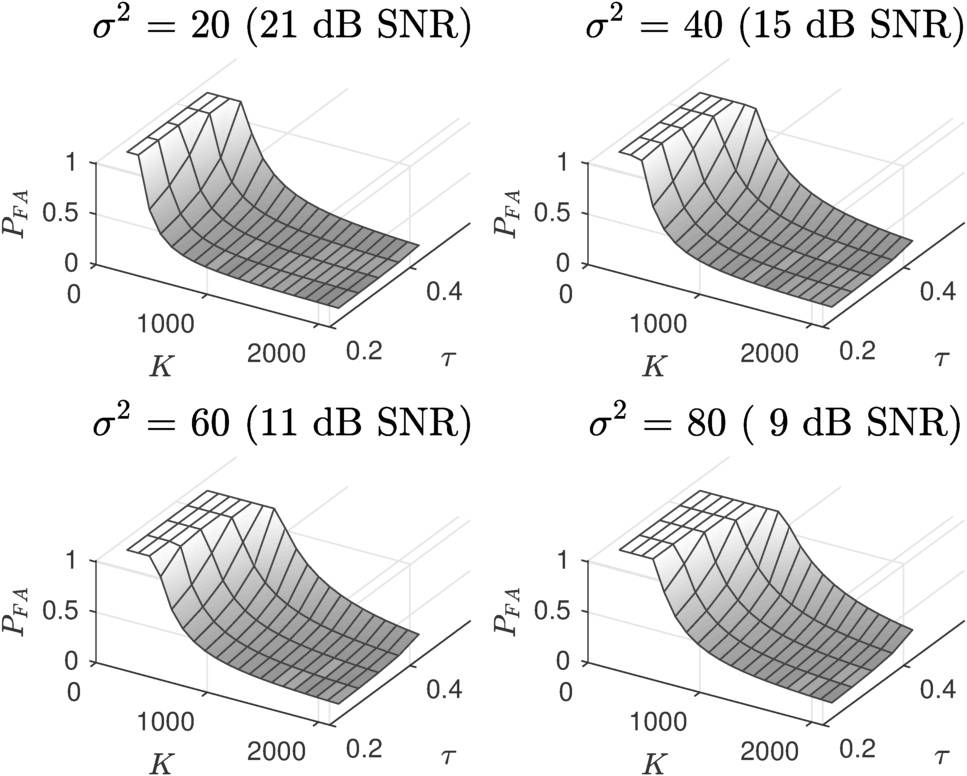}
	\caption{Effect of SNR, $K$, and $\tau$ on the bound on the false alarm probability of $\mathcal{FRDE}$~\eqref{eqn: falseAlarm}.}
	\label{fig: FA1}
\end{figure}

Figure~\ref{fig: FA1} shows how changing the choice of $K$ and $\tau$ affect the upper bound on the false alarm rate of $\mathcal{FRDE}$. We compute the false alarm bound using~\eqref{eqn: falseAlarm} for different values of the local noise covariance $\sigma^2$.
The false alarm probability bound decreases with increasing $K$ and decreasing $\tau$. As we increase the noise covariance $\sigma^2$, we require larger $K$ and/or smaller $\tau$ to achieve the same bound on false alarm probability.

The upper bound on false alarm provided by~\eqref{eqn: falseAlarm} is conservative. Figure~\ref{fig: FA1} shows that, for SNR in the $9$ to $21$ dB range, we require $K$ values ranging from $1000$ to $2000$ to achieve near $0$ upper bound on the false alarm probability. In practice, we can choose $K$ to be much smaller to achieve low false alarm rates. We compute the empirical false alarm rates for different SNR as a function of $K$ and $\tau$. Our simulation considers $4$ different SNR ($9$ dB, 11 dB, 15 dB, 21 dB, corresponding to the SNR considered in Figure~\ref{fig: FA1}). We consider four different values of $\tau$ (0.30, 0.35, 0.40, 0.45), and we consider $K$ in the range from $0$ to $8$. For each level of SNR and each setting of $\tau$ and $K$, we run $100$ simulations. In each simulation, we run $\mathcal{FRDE}$ for $1500$ iterations (with no adversarial agents), and we report the false alarm rate as the ratio of the number of simulations in which any agent reports an adversary to the total number of simulations. 

\begin{figure}[h!]
	\centering
	\includegraphics[keepaspectratio = true, scale = .6]{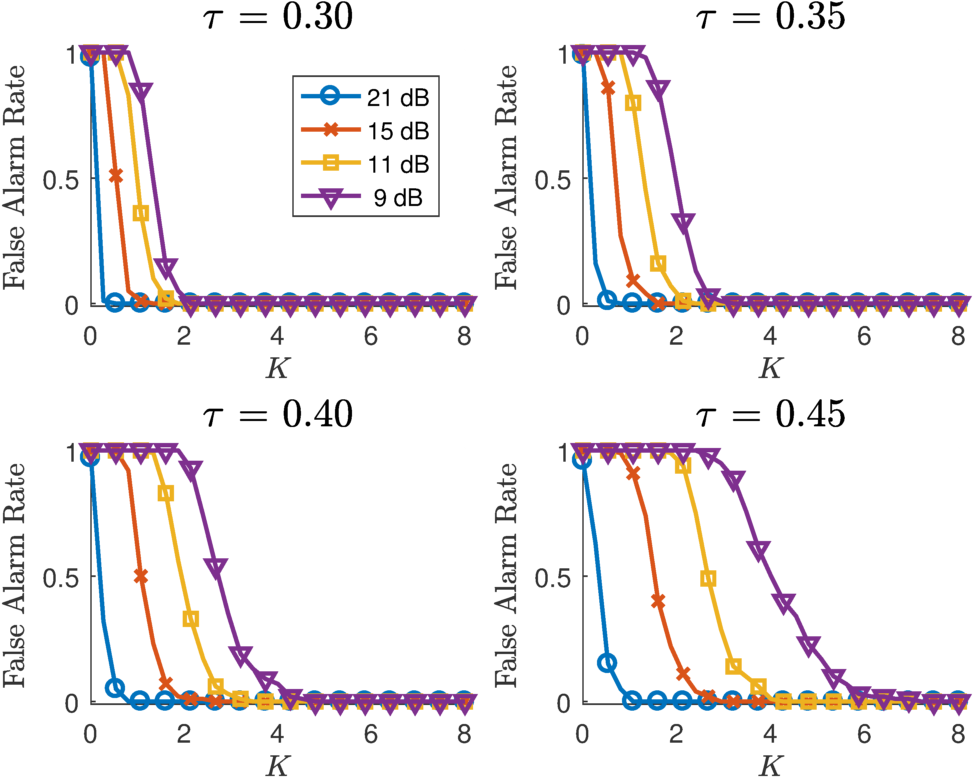}
	\caption{Effect of SNR, $K$, and $\tau$ on empirical false alarm rate of $\mathcal{FRDE}$.}
	\label{fig: FA2}
\end{figure}

Figure~\ref{fig: FA2} shows the effect of SNR, $K$, and $\tau$ on the empirical false alarm rate of $\mathcal{FRDE}$. 
The empirical false alarm rate follows the same trends as the upper bound on false alarm probability from~\eqref{eqn: falseAlarm}. For higher noise covariance (lower SNR), we require larger $K$ and smaller $\tau$ to achieve the same false alarm rate. The difference between empirical rate and the upper bound~\eqref{eqn: falseAlarm}, is that we can achieve low false alarm rates in practice with $K \leq 8$. In contrast, to achieve a low upper bound on false alarm probability, we require $K \geq 1000$. 

\begin{figure}[h!]
	\centering
	\includegraphics[keepaspectratio = true, scale = .6]{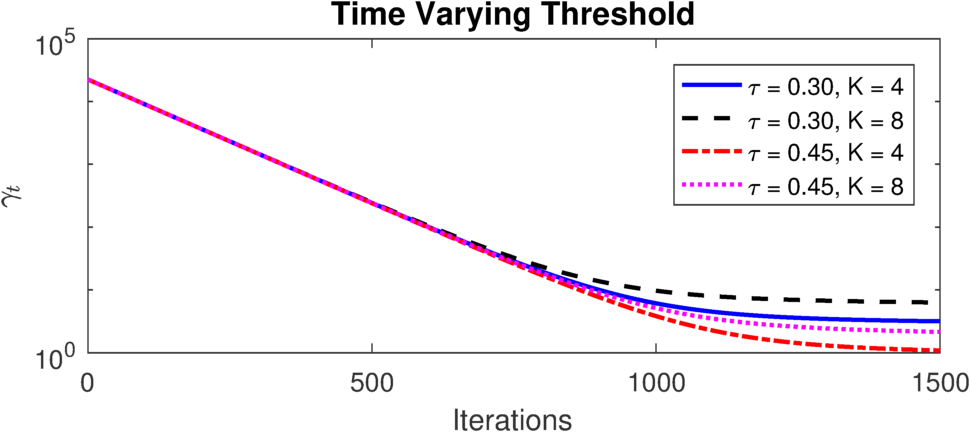}
	\caption{Effect of $K$ and $\tau$ on the bound on the false alarm probability of $\mathcal{FRDE}$.}
	\label{fig: FA3}
\end{figure}

Recall, from~\eqref{eqn: gammaDef}, a smaller value of $K$ results in a smaller value of $\gamma_t$. 
Figure~\ref{fig: FA3} shows how the threshold $\gamma_t$ evolves over time (iterations) for four different choices of $K$ and $\tau$. 
The evolution of $\gamma_t$ does \textit{not} depend on the noise covariance $\sigma^2$. Figure~\ref{fig: FA3} shows that $\gamma_t$ decays more quickly for larger $\tau$, and, for the same value of $\tau$, has smaller value for smaller values of $K$. As a consequence Lemma~\ref{lem: tvSystem}, $\gamma_t$ goes to $0$ eventually (as $t \rightarrow \infty$) for every choice of $K > 0$ and $0 < \tau < \frac{1}{2}$. In practice, however, we are not able to run arbitrarily many iterations of $\mathcal{FRDE}$, so we are interested in the value of $\gamma_t$ after a finite number of iterations.

To achieve small false alarm probabilities, following Figures~\ref{fig: FA1} and~\ref{fig: FA2}, we require smaller $\tau$ and larger $K$. As Figure~\ref{fig: FA3} shows, smaller $\tau$ and larger $K$ yield larger values of $\gamma_t$, which allows adversarial agents to launch more disruptive attacks while evading detection. There is a performance trade-off in~$\mathcal{FRDE}$ between the algorithm's false alarm probability and how disruptive an attack can be before it is detected. 

\section{Conclusion}\label{sect: conclusion}
In this paper, we have studied resilient distributed estimation of a vector parameter by a network of agents. The true but unknown parameter is known to belong to a specified compact subset of a Euclidean space. We have presented an algorithm, Flag Raising Distributed Estimation ($\mathcal{FRDE}$), that allows a network of agents to reliably estimate an unknown parameter in the presence of misbehaving, adversarial agents. Each agent iteratively updates its own estimate based on its previous estimate, its noisy sensor measurement of the parameter, and its neighbors' estimates. An agent raises a flag to indicate the presence of an adversarial agent if any of its neighbors estimates deviates from its own estimate beyond a given threshold. 

Under the global observability condition for the connected normally behaving agents, if the $\mathcal{FRDE}$ algorithm does not detect an attack, then the normally behaving agents correctly estimate the target parameter. The false alarm probability of the algorithm may be arbitrarily small with proper selection of parameters. We have demonstrated the performance of the $\mathcal{FRDE}$ algorithm through numerical examples. Future work includes extending $\mathcal{FRDE}$ to nonlinear measurement models and imperfect communication models.

\appendix
\subsection{Proof of Lemma~\ref{lem: JPD}}
\begin{IEEEproof}
	By inspection, for any $\alpha, \beta > 0$, $J_{\beta, \alpha}$ is the sum of symmetric positive semidefinite matrices, so $J_{\beta, \alpha}$ itself must be symmetric positive semidefinite. To show its positive definiteness, we show that $\mu^T J_{\beta, \alpha} \mu \neq 0$ for all nonzero $\mu \in \mathbb{R}^{NM}$. 

We resort to contradiction. Suppose, there exists a nonzero $\mu \in \mathbb{R}^{NM}$ such that $\mu^T J_{\beta, \alpha} \mu = 0$, which means that $\mu^T \left(L \otimes I_M \right) \mu = 0$ and $\mu^T D_H^T D_H \mu = 0$. Since, by Assumption~\ref{ass: connectivity}, $G$ is connected, $\mu^T \left(L \otimes I_M \right) \mu = 0$ implies that $\mu = 1_N \otimes \overline{\mu}$ for some nonzero $\overline{\mu} \in \mathbb{R}^M$. Then, we have
\begin{equation}\label{eqn: JPDProof1}
	\mu^T D_H^T D_H \mu = \overline{\mu}^T \left(\sum_{n = 1}^N H_n^T H_n \right) \overline{\mu} = 0,
\end{equation}
for some nonzero $\overline{\mu} \in \mathbb{R}^M$. This is a contradiction, since, by Assumption~\ref{ass: globalObservable}, the matrix $\sum_{n=0}^N H_n^T H_n$ is invertible. Thus, we have that the matrix $J_{\beta, \alpha}$ is positive definite. 
\end{IEEEproof}

\subsection{Proof of Lemma~\ref{lem: proc1Conditions}}
\begin{IEEEproof}
	In Procedure 1, the auxiliary parameters $\widehat{\alpha}$ and $\widehat{\beta}$ are chosen to be any positive values. Following Lemma~\ref{lem: JPD}, $J_{\widehat{\beta}, \widehat{\alpha}} \succ 0$, which means that $\lambda_{\max} \left( J_{\widehat{\beta}, \widehat{\alpha}} \right) > 0$. Setting $\alpha$ and $\beta$ according to~\eqref{eqn: 1ab} ensures that $ \lambda_{\max} \left(J_{\beta, \alpha}\right) = 1$ and $0 < \lambda_{\min} \left( J_{\beta, \alpha} \right) \leq 1$. 
\end{IEEEproof}

\subsection{Proof of Lemma~\ref{lem: proc2Conditions}}
\begin{IEEEproof}
	For $\kappa > \frac{\lambda_{\min}(\mathcal{G})}{\lambda_2 (L)}$, let $J^'_{\kappa} = \kappa \left(L \otimes I_M \right) + D_H^T D_H$. We find a lower bound on $\lambda_{\min}\left( J^'_{\kappa} \right)$. The minimum eigenvalue satisfies~\cite{Matrices2}
\begin{equation}\label{eqn: minEig}
	\begin{array}{ccl}\lambda_{\min} \left( J^'_{\kappa}\right) = & \min & z^T J^'_{\kappa} z. \\
	& \text{s.t.}   \left\lVert z \right\rVert = 1 &\end{array}
\end{equation}
Define the subspace
\begin{equation}\label{eqn: consensusSubspace}
	\mathcal{C} = \left\{v \in \mathbb{R}^{NM} \Big\vert v = \mathbf{1}_N \otimes u, u \in \mathbb{R}^M \right\},
\end{equation}
and let $\mathcal{C}^{\perp}$ be the subspace orthogonal to $\mathcal{C}$. The subspace $C$ is in the null space of $L\otimes I_M$. For any, $v \in \mathcal{C}$ with $\left \lVert v \right \rVert = 1$, we can write $v$ as $v = \frac{\mathbf{1}_N \otimes u}{\sqrt{N}}$ for some $u \in \mathbb{R}^M$ with $\left\lVert u \right \rVert = 1$. Then, we can write any $z \in \mathbb{R}^{NM}$ with $\left \lVert z \right \rVert$ as
\begin{equation}\label{eqn: zDecomp}
	z = b \frac{\mathbf{1}_N \otimes u}{\sqrt{N}} + w,
\end{equation}
where $b^2 \leq 1$, $\left\lVert u \right \rVert = 1$, $w \in \mathcal{C}^{\perp}$, and $\left \lVert w \right\rVert^2 = 1- b^2$. 

From~\eqref{eqn: minEig}, we have
\begin{align}
	\begin{split}\label{eqn: bound1}
	& \lambda_{\min} \left( J^'_{\kappa} \right) \geq w^T \left(\kappa\left( L \otimes I_M \right) + D_H^T D_H \right) w + \\
	& \quad \frac{b^2}{N} \left( \mathbf{1}_N \otimes u \right)^T D_H^T D_H \left( \mathbf{1}_N \otimes u \right) + \\
	&  \quad \frac{2b}{\sqrt{N}} \left(\mathbf{1}_N \otimes u \right)^T D_H^T D_H w,
	\end{split}\\
	\begin{split}\label{eqn: bound2}
	& \geq \lambda_2 \left(L\right) \kappa (1 - b^2) - 2 \left\lvert \left \langle w, \frac{b D_H^T D_H \left( \mathbf{1}_N \otimes u \right)}{\sqrt{N}}\right\rangle \right \rvert + \\
	& \quad b^2 \lambda_{\min} \left(\mathcal{G}\right).
	\end{split}
\end{align}
To derive~\eqref{eqn: bound2} from~\eqref{eqn: bound1}, we have used the fact that
\begin{equation}
\begin{split}\label{eqn: bound3}
	\frac{b^2}{N} \left( \mathbf{1}_N \otimes u \right)^T D_H^T D_H \left( \mathbf{1}_N \otimes u  \right) & = b^2 u^T \mathcal{G} u.
\end{split}
\end{equation}
From~\eqref{eqn: bound3} and the Cauchy-Schwarz Inequality, we have
\begin{align}\label{eqn: bound4}
\left\lvert \left \langle w, \frac{b D_H^T D_H \left( \mathbf{1}_N \otimes u \right)}{\sqrt{N}}\right\rangle \right \rvert &\leq \sqrt{\left \lVert w \right \rVert^2b^2 u^T \mathcal{G} u}, \\
	&= \sqrt{\lambda_{\min} \left(\mathcal{G}\right) b^2 \left(1-b^2\right)}. \label{eqn: bound5}
\end{align}
Define the functions $f_1(s)$, $f_2(s)$ as
\begin{align}
	f_1(s) &= \kappa \lambda_2(L) (1-s) + \lambda_{\min}(\mathcal{G}) s, \label{eqn: f1def}\\
	f_2(s) &= f_1(s) - 2 \sqrt{\lambda_{\min}(\mathcal{G}) s(1-s)}. \label{eqn: f2def}
\end{align}
Substituting~\eqref{eqn: bound5} into~\eqref{eqn: bound2}, we have
\begin{equation}\label{eqn: bound6}
	\lambda_{\min} \left( J^'_{\kappa}\right) \geq f_2 \left( b^2 \right).
\end{equation}

We now minimize $f_2(s)$, which has first derivative
\begin{equation}
	\frac{df_2}{ds} = \lambda_{\min}(\mathcal{G})  - \kappa \lambda_2(L) - \frac{\lambda_{\min}(\mathcal{G}) (1 - 2s)}{\sqrt{\lambda_{\min}(\mathcal{G}) s(1-s)}},
\end{equation}
and second derivative
\begin{equation}
	\frac{d^2 f_2}{ds^2} = \frac{\sqrt{\lambda_{\min}(\mathcal{G})}}{2 (s - s^2)}.
\end{equation}
For $0 \leq s \leq 1$, $\frac{d^2 f}{ds^2} \geq 0$, so $f_2(s)$ is convex and minimized for the value of $s$ such that $\frac{df_2}{ds} = 0$ for $0 \leq s \leq 1$. When $\kappa \geq \frac{\lambda_{\min}(G)}{\lambda_2 (L)}$, $\frac{df_2}{ds} = 0$ at
\begin{equation}\label{eqn: minS}
	s =\frac{1}{2} + \frac{1}{2}\sqrt{\frac{\left(\lambda_{\min} \left(\mathcal{G}\right) - \kappa \lambda_2 \left( L \right) \right)^2}{4 \lambda_{\min} \left(\mathcal{ G }\right) + \left(\lambda_{\min} \left(\mathcal{G}\right) - \kappa \lambda_2 \left( L \right) \right)^2 }}.
\end{equation}
By Assumption~\ref{ass: globalObservable}, $\lambda_{\min}(\mathcal{G}) > 0$, so, following~\eqref{eqn: minS}, the minimizing $s$ is less than $1$. Note that, for $0 \leq s < 1$ and $\kappa > \frac{\lambda_{\min} (\mathcal{G})}{\lambda_2(L)}$, we have $f_1(s) > \lambda_{\min} \left(\mathcal{ G} \right)$. Substituting~\eqref{eqn: minS} and $f_1(s) > \lambda_{\min} \left(\mathcal{ G }\right)$ into~\eqref{eqn: f2def} and~\eqref{eqn: bound6}, we have, after algebraic manipulations
\begin{equation}\label{eqn: f2Bound}
\begin{split}
	\lambda_{\min} \left(J^'_{\kappa} \right) &> \lambda_{\min}\left(\mathcal{G}\right)  - \\&\quad \frac{4 \lambda_{\min} \left(\mathcal{G}\right)}{\sqrt{4\lambda_{\min} \left(\mathcal{G}\right) + \left(\lambda_{\min} \left(\mathcal{G}\right) - \lambda_2 \left(L\right) \kappa \right)^2}} 
\end{split}
\end{equation}

Consider the parameters $\alpha$, $\beta$, and $r_1$ selected following Procedure 2. We first show that $\lambda_{\max} \left(J_{\beta, \alpha} \right) \leq 1$. Let $\kappa_1$ be selected according to Step 1) of Procedure 2. Then, we have $\kappa_1 = \frac{\beta}{\alpha}$, and $J_{\beta, \alpha} = \alpha J^'_{\kappa_1}$. The maximum eigenvalue of of $J^'_{\kappa_1}$ satisfies
\begin{equation}\label{eqn: tempMaxBound}
	\lambda_{\max} \left( J^'_{\kappa_1}\right) \leq \kappa_1 \lambda_{\max} \left( L \right) + \lambda_{\max} \left( D_H^T D_H \right).
\end{equation}
As a consequence of Assumption~\ref{ass: measureMatrix}, $\lambda_{\max} \left( D_H^T D_H \right) \leq N$. From~\eqref{eqn: tempMaxBound}, we see that choosing $\alpha$ according to Procedure 2 ensures that $\lambda_{\max} \left( J_{\beta, \alpha} \right) \leq 1$. 

We now show that $0 < r_1 \leq \lambda_{\min} \left( J_{\beta, \alpha} \right)$. As a consequence of Assumption~\ref{ass: measureMatrix}, we have $\lambda_{\min} (\mathcal{G}) \leq 1$. By algebraic manipulation, the selection of $\kappa_1$ in Step 1) ($\kappa_1 > \frac{1}{\lambda_2 (L)} \left(\lambda_{\min} (\mathcal{G}) + 2 \sqrt{4 - \lambda_{\min} (\mathcal{G})}\right)$) ensures that $r_1 > 0$. Since $J_{\beta, \alpha} = \alpha J^'_{\kappa_1}$, we have $\lambda_{\min} \left( J_{\beta, \alpha} \right) = \alpha \lambda_{\min} \left( J^'_{\kappa_1}\right)$. Note that the right hand side of~\eqref{eqn: f2Bound} is equal to $\frac{r_1}{\alpha}$. Thus, from~\eqref{eqn: f2Bound}, we have $r_1 \leq \lambda_{\min} \left( J_{\beta, \alpha} \right)$. 
\end{IEEEproof}

\subsection{Proof of Lemma~\ref{lem: timeAveragedNoise}}
\begin{IEEEproof}
	Define the process $\left\{V_t\right\}$ as
\begin{equation}\label{eqn: vtDef}
	\begin{split}
		V_{t} &= \left(\frac{t}{t+1} \right)^{1-\tau_0} V_{t-1} + \frac{1}{(t+1)^{1-\tau_0}} w_t, \\
		V_0 &= m_0 = w_0.
	\end{split}
\end{equation}
for $0 \leq \tau_0 < {1 \over 2}$. Note that $V_t = (t+1)^{\tau_0} m_t$. Define the process $\left\{\widetilde{V}_t\right\}$ as
\begin{equation}\label{eqn: vtTildeDef}
	\widetilde{V}_t = \left\lVert V_t \right\rVert^2 +  \sum_{j = t+1}^{\infty} {\trace\left(\Sigma\right) \over (j+1)^{2(1-\tau_0)}}.
\end{equation}
By definition, $\widetilde{V}_t \geq 0$. We now show that $\left\{ \widetilde{V}_t \right\}$ is a supermartingale. 

Substituting~\eqref{eqn: vtDef} into~\eqref{eqn: vtTildeDef} and performing algebraic manipulations, we have
\begin{equation}\label{eqn: taProof1}
	\begin{split}
		\widetilde{V}_{t+1} &= \frac{\left\lVert w_{t+1} \right\rVert^2 + 2 w_{t+1}^T V_t (t+1)^{1-\tau_0}}{(t+2)^{2(1-\tau_0)}} + \\
			&\left\lVert V_t \right\rVert^2\left(\frac{t+1}{t+2}\right)^{2(1-\tau_0)} + \sum_{j = t+2}^{\infty} {\trace\left(\Sigma\right) \over (j+1)^{2(1-\tau_0)}}.
	\end{split}
\end{equation}
The processes $V_t$ and $\widetilde{V}_t$ depend only on $w_0, \dots, w_t$, which means that $w_{t+1}$ is independent of $V_0, \dots, V_t$ and $\widetilde{V}_0, \dots, \widetilde{V}_t$. From~\eqref{eqn: vtTildeDef}, we also have
\begin{equation}\label{eqn: taProof2}
	\mathbb{E}\left[\left\lVert V_t\right\rVert^2 \Big\vert \widetilde{V}_0, \dots, \widetilde{V}_t\right] = \left\lVert V_t\right\rVert^2. 
\end{equation}
Taking the expectation of~\eqref{eqn: taProof1} conditioned on $\widetilde{V}_0, \dots, \widetilde{V}_t$, we have
\begin{align}
	\begin{split}\label{eqn: taProof3}
		\mathbb{E}& \left[\widetilde{V}_{t+1} \Big\vert \widetilde{V}_0, \dots, \widetilde{V}_t \right] = \frac{\trace\left(\Sigma\right)}{(t+2)^{2(1-\tau_0)}} + \\ &\quad \left\lVert V_t \right\rVert^2\left(\frac{t+1}{t+2}\right)^{2(1-\tau_0)} + \sum_{j = t+2}^{\infty} {\trace\left(\Sigma\right) \over (j+1)^{2(1-\tau_0)}},
	\end{split} \\
	&\leq \left\lVert V_t \right\rVert^2 + \sum_{j = t+1}^{\infty} {\trace\left(\Sigma\right) \over (j+1)^{2(1-\tau_0)}} = \widetilde{V}_t, \label{eqn: taProof4}
\end{align}
where~\eqref{eqn: taProof4} follows from~\eqref{eqn: taProof3} since $\left(\frac{t+1}{t+2} \right)^{2(1-\tau_0)} \leq 1$. 

First, we prove~\eqref{eqn: averageNoiseConvergence}. Since $\left\{ \widetilde{V}_t \right\}$ is a nonnegative supermartingale, it converges almost surely to a finite, nonnegative random variable $V^*$, i.e.
\begin{equation}\label{eqn: taProof10}
	\lim_{t \rightarrow \infty} \widetilde{V}_t = V^* \text{ a.s.}.
\end{equation}
Since $\lim_{t \rightarrow \infty} \sum_{j = t+1}^{\infty} {\trace\left(\Sigma\right) \over (j+1)^{2(1-\tau_0)}} = 0$, from~\eqref{eqn: vtTildeDef}, we also have
\begin{equation}\label{eqn: taProof11}
	\lim_{t \rightarrow \infty} \left\lVert V_t \right\rVert^2 = V^* \text{ a.s.}.
\end{equation}
For a finite, nonnegative random variable $V*$, we have $V^* = 0 \text{ a.s.}$ if and only if $\mathbb{E} \left[ V^* \right] = 0$. By Fatou's Lemma, we have
\begin{equation}\label{eqn: taProof12}
	\mathbb{E} \left[ V^* \right] \leq \liminf_{t \rightarrow \infty} \mathbb{E} \left[ \left\lVert V_t \right \rVert^2 \right].
\end{equation}

To evaluate the right hand side of~\eqref{eqn: taProof12}, note that $\left\lVert V_t \right \rVert^2 = \left(1 + t \right)^{\tau_0} \left \lVert m_t \right \rVert^2$. From~\eqref{eqn: averageNoise}, we can express $\mathbb{E} \left[ \left\lVert m_{t+1} \right  \rVert^2 \right]$ as
\begin{align}\label{eqn: taProof13}
	\mathbb{E}& \left[ \left\lVert m_{t+1} \right \rVert^2 \right] = \mathbb{E} \left[ \left\lVert\left(1- \frac{1}{t+1} \right)m_{t} + \frac{w_{t+1}}{t+1} \right \rVert^2\right], \\
	\begin{split}\label{eqn: taProof14}
		=& \left(1 - \frac{1}{t+1}\right)^2 \mathbb{E} \left[ \left\lVert m_t \right \rVert^2 \right] + \frac{1}{(t+1)^2} \mathbb{E} \left[\left\lVert w_{t+1} \right\rVert^2\right],
	\end{split}\\
	\begin{split}\label{eqn: taProof15}
		=& \left(\!1 - \!\left({2 \over t+1} - {1\over (t+1)^2} \right)\right)\!\mathbb{E} \left[ \left\lVert m_t \right \rVert^2 \right] \!+ \!{\trace \left(\Sigma \right) \over (t+1)^2}\!,
	\end{split}\\
	\leq & \left(1 - \frac{1}{t+1} \right) \mathbb{E} \left[ \left \lVert m_t \right\rVert^2 \right] + {\trace\left(\Sigma\right) \over (t+1)^2}, \label{eqn: taProof16}
\end{align}
where~\eqref{eqn: taProof14} follows from~\eqref{eqn: taProof13} since $w_{t+1}$ is independent of $m_t$, and~\eqref{eqn: taProof16} follows from~\eqref{eqn: taProof15} since $(t+1)^{-1} \leq 2(t+1)^{-1} - (t+1)^{-2}$ for all $t \geq 0$. Relation~\eqref{eqn: taProof16} falls under the purview of Lemma~\ref{lem: tvSystem}, which means that 
\begin{align}\label{eqn: taProof17}
	\lim_{t \rightarrow \infty} (t+1)^{2\tau_0} \mathbb{E} \left[ \left \lVert m_t \right\rVert^2 \right] = 0.
\end{align}
Since $\left\lVert V_t \right\rVert = (1+t)^{\tau_0} \left\lVert m_t \right\rVert$, and since $V^* \geq 0$, combining~\eqref{eqn: taProof11} with~\eqref{eqn: taProof17} yields
\begin{equation}\label{eqn: taProof18}
	\mathbb{E} \left[ V^* \right] = 0,
\end{equation}
which means that $V^* = 0$ a.s.. That is, $\left\lVert V_t \right \rVert^2$, converges to $0$ almost surely.  Substituting $\left\lVert V_t \right \rVert = (1+t)^{\tau_0} \left\lVert m_t \right \rVert$, we have
\begin{equation}\label{eqn: taProof19}
	\mathbb{P} \left( \lim_{t \rightarrow \infty} (t+1)^{\tau_0} \left \lVert m_t \right\rVert = 0 \right) = 1,
\end{equation}
for every $0 \leq \tau_0 < \frac{1}{2}$, which shows~\eqref{eqn: averageNoiseConvergence}. 

Second, we prove~\eqref{eqn: averageNoiseSup}. By~\eqref{eqn: taProof4}, we have that $\left\{\widetilde{V}_t\right\}$ is a nonnegative supermartingale. Then, by the maximal inequality for nonnegative supermartingales~\cite{Kushner}, we have, for any $k > 0$,
\begin{equation}\label{eqn: taProof5}
	\mathbb{P} \left(\sup_{t \geq 0} \widetilde{V}_t > k^2 \right) \leq \frac{1}{k^2}\mathbb{E} \left[ \widetilde{V}_0 \right].
\end{equation}
From~\eqref{eqn: vtDef} and~\eqref{eqn: vtTildeDef}, we have
\begin{align}\label{eqn: taProof6}
	\mathbb{E} \left[ \widetilde{V}_0 \right] &= \mathbb{E} \left[\left\lVert V_0 \right\rVert^2 \right] + \sum_{j = 2}^{\infty} {\trace\left(\Sigma\right) \over j^{2(1-\tau_0)}}, \\
	& = \sum_{j = 1}^{\infty} {\trace\left(\Sigma\right) \over j^{2(1-\tau_0)}},\label{eqn: taProof7}
\end{align}
where~\eqref{eqn: taProof7} follows from~\eqref{eqn: taProof6} since $V_0 = m_0 = w_0$ and $\mathbb{E} \left[ \left\lVert w_0 \right\rVert^2 \right] = \trace\left(\Sigma\right)$. 

From~\eqref{eqn: vtTildeDef}, we have $\left\lVert V_t \right \rVert^2 \leq \widetilde{V}_t$, and, from~\eqref{eqn: vtDef}, we have $\left\lVert V_t \right \rVert^2 = (1 + t)^{2\tau_0} \left\lVert m_t \right \rVert^2$. Thus, we have
\begin{equation}\label{eqn: taProof8}
	\mathbb{P} \left(\sup_{t \geq 0} \widetilde{V}_t > k^2 \right) \geq \mathbb{P} \left(\sup_{t\geq 0} (t+1)^{2\tau_0} \left\lVert m_t \right\rVert^2 > k^2 \right).
\end{equation}
Combining~\eqref{eqn: taProof5},~\eqref{eqn: taProof7}, and~\eqref{eqn: taProof8}, we have, after algebraic manipulations
\begin{equation}\label{eqn: taProof9}
	\mathbb{P} \left(\sup_{t \geq 0 } \left \lVert m_t \right \rVert^2 > \frac{k}{(t+1)^{\tau_0}}\right) \leq {1 \over k^2} \sum_{j = 1}^{\infty} {\trace\left(\Sigma\right) \over j^{2(1-\tau_0)}},
\end{equation}
which establishes~\eqref{eqn: averageNoiseSup}.

\end{IEEEproof}

\subsection{Proof of Lemma~\ref{lem: LPartition}} 
\begin{IEEEproof}
	First, we show that $J^{\mathcal{X}}_{\beta, \alpha}$ is positive definite. Since $\beta, \alpha > 0$, the graph $G_{\mathcal{X}}$ is connected, and $\mathcal{X}$ is globally observable, Lemma~\ref{lem: JPD} applies, and, we have that the matrix $J^{\mathcal{X}}_{\beta, \alpha}$ is positive definite. 

Second, we show that $\lambda_{\max} \left( J^{\mathcal{X}}_{\beta, \alpha} \right) \leq 1$. Without loss of generality, let $\mathcal{X} = \left\{1, \dots, \left\vert \mathcal{X} \right\vert \right\},$ and let $\mathcal{Y} = V\setminus \mathcal{X}$. Then, we can partition the matrix $J_{\beta, \alpha}$ as
\begin{equation}\label{eqn: PartitionProof1}
	J_{\beta, \alpha} = \left[\begin{array}{cc} \widetilde{J}^{\mathcal{X}}_{\beta, \alpha} & \left( K^1_{\beta, \alpha} \right)^T \\ K^1_{\beta, \alpha} & K^2_{\beta, \alpha} \end{array} \right],
\end{equation}
where 
\begin{equation}\label{eqn: tildeJDef}
	\widetilde{J}^{\mathcal{X}}_{\beta, \alpha} = J^{\mathcal{X}}_{\beta, \alpha} + \beta\left(\Sigma_{\mathcal{X}, \mathcal{Y}} \otimes I_M\right).
\end{equation}
We show that $\lambda_{\max} \left( \widetilde{J}^{\mathcal{X}}_{\beta, \alpha} \right) \leq 1$. For purposes of contradiction, suppose that $\lambda_{\max} \left( \widetilde{J}^{\mathcal{X}}_{\beta, \alpha} \right) > 1$, and let $v$ be the associated eigenvector. Thus, we have
\begin{equation}\label{eqn: PartitionProof2}
	\sqrt{ v^T \left.\widetilde{J}^{\mathcal{X}}_{\beta, \alpha}\right.^T \left. \widetilde{J}^{\mathcal{X}}_{\beta, \alpha}\right. v} = \lambda_{\max} \left( \widetilde{J}^{\mathcal{X}}_{\beta, \alpha}\right) \left\lVert v \right\rVert > \left\lVert v \right \rVert.
\end{equation}
By definition, we have
\begin{equation}\label{eqn: PartitionProof3}
	\lambda_{\max}\left(\widetilde{J}_{\beta, \alpha} \right) = \sup_{\overline{v} \in \mathbb{R}^{NM}} \frac{\left\lVert \widetilde{J}_{\beta, \alpha} \overline{v} \right \rVert}{\left\lVert \overline{v} \right\rVert}.
\end{equation}
Consider $\widetilde{v} = \left[\begin{array}{cc} v^T & 0_{\left(N -\left\vert \mathcal{X} \right\vert \right) M}^T \end{array} \right]^T$. For the vector $\widetilde{v} \in \mathbb{R}^{NM}$, we have
\begin{align}
	\left\lVert {J}_{\beta, \alpha} \widetilde{v} \right\rVert &= \sqrt{{v}^T  \left( \left.\widetilde{J}^{\mathcal{X}}_{\beta, \alpha}\right.^T \widetilde{J}^{\mathcal{X}}_{\beta, \alpha} + {K^1_{\beta, \alpha}}^T K^1_{\beta, \alpha}\right){v}}, \label{eqn: PartitionProof4} \\
		& \geq \sqrt{{v}^T  \left( \left. \widetilde{J}^{\mathcal{X}}_{\beta, \alpha}\right.^T \widetilde{J}^{\mathcal{X}}_{\beta, \alpha}\right){v}},\label{eqn: PartitionProof5}
\end{align}
where~\eqref{eqn: PartitionProof5} follows from~\eqref{eqn: PartitionProof4} because the matrix $ {K^1_{\beta, \alpha}}^T K^1_{\beta, \alpha}$ is positive semidefinite. Then, substituting for~\eqref{eqn: PartitionProof2}, we have
\begin{align}\label{eqn: PartitionProof6}
	\frac{\left\lVert J_{\beta, \alpha} \widetilde{v} \right\rVert}{\left \lVert \widetilde{v} \right\rVert} \geq \frac{\left\lVert \widetilde{J}^{\mathcal{X}}_{\beta, \alpha} v \right\rVert}{\left \lVert v \right\rVert} = \lambda_{\max} \left( \widetilde{J}^{\mathcal{X}}_{\beta, \alpha}\right) > 1.
\end{align}
From~\eqref{eqn: PartitionProof3} and~\eqref{eqn: PartitionProof6}, we have $\lambda_{\max} \left( J_{\beta, \alpha} \right) > 1$, which is a contradiction since $\alpha, \beta$ are chosen such that $\lambda_{\max} \left( J_{\beta, \alpha} \right) \leq 1$. Thus, we have $\lambda_{\max} \left( \widetilde{J}^{\mathcal{X}}_{\beta, \alpha} \right) \leq 1$ as well. Applying Weyl's Inequality~\cite{Matrices} to~\eqref{eqn: tildeJDef}, we have
\begin{equation}\label{eqn: PartitionProof7}
	\lambda_{\max} \left( J^{\mathcal{X}}_{\beta, \alpha}\right) + \lambda_{\min} \left( \beta \left(\Sigma_{\mathcal{X}, \mathcal{Y}} \otimes I_M \right) \right) \leq \lambda_{\max} \left( \widetilde{J}^{\mathcal{X}}_{\beta, \alpha} \right).
\end{equation}
By definition, $\Sigma_{\mathcal{X}, \mathcal{Y}} \succeq 0$, so $\lambda_{\min} \left( \beta \left( \Sigma_{\mathcal{X}, \mathcal{Y}} \otimes I_M \right) \right) \geq 0$. Then, from~\eqref{eqn: PartitionProof7}, we have $\lambda_{\max} \left( J_{\beta, \alpha}^{\mathcal{X}} \right) \leq 1$. 
\end{IEEEproof}

\bibliography{IEEEabrv,References}

\end{document}